\documentclass[11pt,a4paper]{amsart}
\usepackage{amsmath,amssymb,amsthm,enumitem}
\usepackage[a4paper]{geometry}
\usepackage{palatino}
\usepackage{graphicx}
\usepackage{xcolor}
\usepackage[colorlinks=true]{hyperref}
\hypersetup{urlcolor=blue, citecolor=red, linkcolor=blue}
\usepackage{cite}

\numberwithin{equation}{section}

\newcommand{\e}{\varepsilon}
\newcommand{\R}{\mathbb{R}}
\newcommand{\dd}{\,\mathrm{d}}
\newcommand{\bd}{\boldsymbol}
\newcommand{\Zs}{\mathcal Z}
\DeclareMathOperator{\sign}{sign}
\DeclareMathOperator{\meas}{meas}
\DeclareMathOperator{\interior}{int}

\newtheorem{theorem}{Theorem}[section]
\newtheorem{lemma}[theorem]{Lemma}
\newtheorem{definition}[theorem]{Definition}
\newtheorem{proposition}[theorem]{Proposition}
\newtheorem{corollary}[theorem]{Corollary}
\newtheorem{conjecture}[theorem]{Conjecture}
\theoremstyle{definition}
\newtheorem{remark}[theorem]{Remark}

\begin{document}

\title[One bubble against a ring]{One bubble against a ring:\\ new sign-changing solutions of the Brezis-Nirenberg problem}

\author[A. Pistoia]{Angela Pistoia}
\address{\noindent  Dipartimento di Scienze di Base e Applicate per l'Ingegneria, Università degli Studi di Roma Sapienza}
\email{angela.pistoia@uniroma1.it}

\author[G. M. Rago]{Giuseppe Mario Rago}
\address{\noindent  Dipartimento di Matematica, Universit\`a degli Studi di Bari Aldo Moro, Italy }
\email{g.rago6@phd.uniba.it}

\author[G. Vaira]{Giusi Vaira}
\address{\noindent  Dipartimento di Matematica, Universit\`a degli Studi di Bari Aldo Moro, Italy }
\email{giusi.vaira@uniba.it}

\thanks{G.~M. Rago and G.~Vaira acknowledge support from the INdAM-GNAMPA
project ``Singolarità, nonlinearità e interazioni non locali: nuovi
approcci analitici per equazioni differenziali e sistemi''
(CUP E53C25002010001).}


\begin{abstract}
We consider the Brezis--Nirenberg problem $-\Delta u=|u|^{2^*-2}u+\varepsilon u$ in the unit ball of $\mathbb R^N$, $N\ge5$, with Dirichlet boundary conditions, and we look for sign-changing solutions which, as $\varepsilon\to0^+$, blow up positively at the centre and negatively at the $k$ vertices of a regular polygon. The key observation is that, for this configuration, the gradient of the reduced energy factorises, so that the whole finite-dimensional problem is governed by the zeros of a single explicit function of the radius of the polygon. Its behaviour depends on the dimension: when $N=5$ such solutions exist for every $k\ge2$, and the ring of negative peaks approaches the boundary at distance $O(1/k)$; when $N\ge6$ there are two solutions for $2\le k\le6$, while no solution of this type exists for $k$ large. For $2\le k\le6$ the nodal set of every such solution reaches the boundary of the ball, in contrast with what happens for the almost critical problem. Nevertheless, the two solutions found for $N\ge6$ have very different nodal structures: numerically, for $3\le k\le6$ one of them has only two nodal domains, and its nodal set meets the plane of the peaks along a closed curve around the centre, while the other one can have as many as $k+1$ nodal domains. These results extend to higher dimensions a configuration previously found in dimension four.
\end{abstract}

\maketitle

\section{Introduction and main results}

In this paper we continue the study started in \cite{PRV2} of sign-changing blowing-up solutions of the Brezis--Nirenberg problem
\begin{equation}\label{pb}
  \begin{cases}
    -\Delta u=|u|^{2^*-2}u+\e u & \text{in }\Omega,\\
    u=0 & \text{on }\partial\Omega,
  \end{cases}
\end{equation}
where $\Omega\subset\R^N$, $N\ge 5$, is a smooth bounded domain,
$\e>0$ is a small parameter and, as usual, $2^*=\frac{2N}{N-2}$ is the critical Sobolev exponent.
Problem \eqref{pb} was introduced by Brezis and Nirenberg in the celebrated paper \cite{BN}, where it is shown that, if $N\ge4$, it has a positive solution for every $\e\in(0,\lambda_1(\Omega))$, $\lambda_1(\Omega)$ being the first eigenvalue of $-\Delta$ in $\Omega$ with Dirichlet boundary conditions, while in dimension $N=3$ positive solutions exist only for $\e$ in a proper subinterval of $(0,\lambda_1(\Omega))$. On the other hand, by the Pohozaev identity, \eqref{pb} has no nontrivial solutions for $\e\le0$ when $\Omega$ is star-shaped. The lack of compactness of the embedding $H^1_0(\Omega)\hookrightarrow L^{2^*}(\Omega)$ is reflected in the behaviour of the solutions as $\e\to0^+$: they may blow up, and their profile near each blow-up point is described by a \emph{bubble}, i.e.\ an extremal of the Sobolev inequality \cite{A,T}; the bubbles are the only positive solutions of the limit problem in $\R^N$ \cite{CGS}.

For positive solutions this phenomenon is by now well understood. Han \cite{H} and Rey \cite{R} proved that minimal energy solutions blow up, as $\e\to0$, at a critical point of the Robin function of $\Omega$, and Rey \cite{R} showed that, conversely, in dimension $N\ge5$ non-degenerate critical points of the Robin function generate families of solutions blowing up at them. Positive solutions concentrating at several points were constructed by Musso and Pistoia \cite{MP}: in this case the role of the Robin function is played by a finite-dimensional function which involves both the Robin function and the Green function of $\Omega$, and the geometry of the domain becomes crucial. A fine description of multi-bubble blow-up has been obtained by K\"onig and Laurain \cite{KL1,KL2} (see also Druet \cite{D} for the three-dimensional case), and the construction of bubbling solutions in general bounded domains in low dimensions has been recently addressed in \cite{LVWW,LVWW1}.

The picture for sign-changing solutions is much richer and still far from complete. Sign-changing solutions may blow up at different points with different signs \cite{BMP,MiPi,IV2}, or at the same point, as a \emph{tower} of bubbles with alternating signs \cite{IV1,Pre}; nodal radial solutions in the ball, for which the dependence on the dimension is particularly evident, have been studied in \cite{AGGPV,PV}, and related phenomena for critical problems with a Hardy potential in \cite{EGPV}. In all these results the dimension plays a crucial role, since the concentration rate of the bubbles, and hence the structure of the finite-dimensional problem obtained by the reduction, depends on $N$: the scales are of order $\e^{\frac1{N-4}}$ when $N\ge5$, while they are exponentially small in $1/\e$ when $N=4$.

\medskip

\noindent In dimension $N\ge 5$, sign-changing multi-bubble solutions of critical and supercritical problems were studied by Micheletti and Pistoia \cite{MiPi}. For the Brezis--Nirenberg problem their reduction leads to a finite-dimensional function which depends simultaneously on the concentration points and on the individual scales. To describe it, and to explain the difficulties we have to face, let us introduce some notation which will be used throughout the paper.
Let $G(x, y)$ be the Green function of $-\Delta$ in $\Omega$ with Dirichlet boundary conditions and let $H$ be its regular part, i.e.
\[G(x, y)=\frac{1}{\omega_N}\frac{1}{|x-y|^{N-2}}-H(x, y),\]
where $\omega_N:=(N-2)\,\meas(\mathbb S^{N-1})$, $\mathbb S^{N-1}$ is the unit sphere of $\R^N$ and, for every $y\in\Omega$,
\[-\Delta_x H(x, y)=0\quad\text{in}\quad \Omega,\qquad H(x, y)=\frac{1}{\omega_N|x-y|^{N-2}}\quad\text{on}\quad \partial\Omega.\]
The leading term of the regular part, $\tau(x):=H(x,x)$, is called the \emph{Robin function} of $\Omega$.
Given an integer $k\ge 2$, we set
\[\Omega^*_k:=\left\{\bd\xi=(\xi_1, \ldots, \xi_k)\in \Omega^k\,:\, \xi_i\neq \xi_j \quad\text{for}\quad i\neq j\right\}\]
and, for $\bd\xi\in\Omega^*_k$ and $\bd\sigma=(\sigma_1,\ldots,\sigma_k)\in\{-1,1\}^k$, we define the symmetric matrix
\[
M^{\bd\sigma}(\bd\xi)=
\begin{pmatrix}
\tau(\xi_1)&-\sigma_1\sigma_2G(\xi_1,\xi_2)&\cdots&-\sigma_1\sigma_kG(\xi_1,\xi_k)\\
-\sigma_1\sigma_2G(\xi_1,\xi_2)&\tau(\xi_2)&\cdots&-\sigma_2\sigma_kG(\xi_2,\xi_k)\\
\vdots&\vdots&\ddots&\vdots\\
-\sigma_1\sigma_kG(\xi_1,\xi_k)&-\sigma_2\sigma_kG(\xi_2,\xi_k)&\cdots&\tau(\xi_k)
\end{pmatrix}.\]
$M^{\bd\sigma}(\bd\xi)$ is the Green--Robin interaction matrix: its diagonal entries describe the interaction of each bubble with the boundary, while the off-diagonal entries describe the interaction between distinct bubbles.
In our notation, the leading part of the reduced energy in \cite{MiPi} is
\begin{equation}\label{eq:FMiPi}
F^{\bd\sigma}_k(\bd d,\bd\xi)
=
\frac{a_N}{2}\left[\sum_{j=1}^k \tau(\xi_j) d_j^{N-2}-\sum_{\substack{i, j=1\\ i\neq j}}^k \sigma_i\sigma_j G(\xi_i, \xi_j)d_i^{\frac{N-2}{2}}d_j^{\frac{N-2}{2}}\right]-b_N\sum_{j=1}^k d_j^2,
\end{equation}
or equivalently, setting $\kappa_j:=d_j^{\frac{N-2}{2}}$,
\begin{equation}\label{eq:FMiPi1}
F^{\bd\sigma}_k(\bd \kappa,\bd\xi)=\frac{a_N}{2}\left\langle
M^{\bd\sigma}(\bd\xi)\bd \kappa,
\bd \kappa
\right\rangle
-
b_N\sum_{j=1}^k\kappa_j^{\frac{4}{N-2}},
\end{equation}
where $a_N, b_N>0$ are given in \eqref{eq:aNbN}.
Thus the geometry of the domain and the nodal pattern enter the reduced problem simultaneously.

The general reduction allows one to prescribe arbitrary signs, but the analysis of the critical points of \eqref{eq:FMiPi} is rather delicate, because the location variables and all the relative scales remain coupled. This is one of the reasons why explicit constructions in natural geometries have been obtained only for particular configurations. For example, in \cite{MiPi} a family with an even number of bubbles placed at the vertices of a regular polygon with alternating signs is constructed in the unit ball: the symmetry reduces all the scales to one parameter and all the locations to one radial variable. The same paper also produces many sign-changing patterns in suitable dumbbell-shaped domains, and treats the case of one positive and one negative bubble in an arbitrary smooth bounded domain. In all these cases the difficulty is not the reduction itself, but the analysis of the critical points of the finite-dimensional function.
In dimension $N=4$ the situation is different: the concentration scales are exponentially small in $1/\e$ and the leading equations for the relative scales become spectral. Hence the scale variables can be eliminated through the spectral structure of the interaction matrix, which makes accessible configurations that are considerably harder to detect through the higher-dimensional reduced function.

\medskip

This spectral mechanism was exploited in \cite{PRV2}, where several families of sign-changing multi-peak solutions were constructed in dimension four, in configurations whose sign pattern is not alternating and which are therefore not covered by the symmetric constructions of \cite{MiPi}. Among them: in the unit ball, a positive bubble at the centre surrounded by $k$ negative bubbles placed at the vertices of a regular polygon, and positive and negative bubbles placed at the vertices of two regular polygons lying in orthogonal planes; in an annulus, bubbles with alternating signs at the vertices of a regular polygon. Whether these configurations survive in higher dimensions was left open in \cite{PRV2}. When $N\ge5$ the scales are only polynomially small in $\e$, the self-interaction term $b_N\sum_j d_j^2$ in \eqref{eq:FMiPi} is of the same order as the interaction terms, and the scales can no longer be eliminated by a spectral argument: the location of the bubbles and their relative scales must be determined together, as critical points of the full function \eqref{eq:FMiPi}.

The aim of this paper is to show the persistence of the multi-peak configurations constructed in \cite{PRV2} for $N=4$ in dimensions $N\ge 5$, where no spectral mechanism is available and one has to deal with the full function \eqref{eq:FMiPi}. We focus on a particularly interesting one: in the unit ball we look for solutions with one positive peak at the centre, surrounded by an arbitrary number $k\ge2$ of negative peaks placed at the vertices of a regular polygon.

To state our results we introduce the relevant symmetries. Let $\Gamma_k\subset O(N)$ be the group generated by the rotation of angle $2\pi/k$ in the $(x_1,x_2)$-plane, by the reflection $(x_1,x_2,x'')\mapsto(x_1,-x_2,x'')$ and by $O(N-2)$ acting on $x''=(x_3,\dots,x_N)$, and let
\[
\mathcal H_k:=\{u\in H^1_0(\mathcal B_1):\ u\circ\gamma=u\quad\forall\gamma\in\Gamma_k\},
\]
where $\mathcal B_1$ denotes the unit ball of $\R^N$ centred at the origin.
For $r\in(0,1)$ we consider the vertices of a regular $k$-gon of radius $r$ in the $(x_1,x_2)$-plane,
\begin{equation}\label{eq:xi}
\xi_j=\xi_j(r):=r\Big(\cos\tfrac{2\pi(j-1)}k,\ \sin\tfrac{2\pi(j-1)}k,\ 0,\dots,0\Big),\qquad j=1,\dots,k,
\end{equation}
and we set $\mathcal M:=(0,1)\times(0,\infty)^2$. We look for solutions of the form
\begin{equation}\label{eq:form}
u_\e=PU_{\delta_0,0}-\sum_{j=1}^kPU_{\delta_1,\xi_j(r)}+\phi_\e,\qquad \delta_i=d_i\,\e^{\frac1{N-4}},\quad i=0,1,
\end{equation}
where $PU_{\delta,\xi}$ is the projection on $H^1_0(\mathcal B_1)$ of the standard bubble (see Section~\ref{sec:prelim}), $(r,d_0,d_1)=(r_\e,d_{0,\e},d_{1,\e})\in\mathcal M$ and $\|\phi_\e\|_{H^1_0}\to0$ as $\e\to0$.

The key observation of the paper (Lemma~\ref{lem:factor}) is that, for this configuration, the gradient of the reduced function factorises as the product of an invertible matrix and of a vector whose components can be solved one at a time. As a result, the critical points of the reduced function, and their local Brouwer degree, are completely described by the zeros and the sign changes of a single explicit function $g:(0,1)\to\R$ of the radius of the polygon, defined in \eqref{eq:g}. This function always tends to $+\infty$ as $r\to0^+$, while as $r\to1^-$ it tends to $-\infty$ if $N=5$ and to $1$ if $N\ge6$ (Lemmas~\ref{lem:zero} and~\ref{lem:one}). This dichotomy, which is due to the sign of the exponent $\mathfrak e_N=\frac{N^2-6N+4}{N-2}$, is at the origin of the different results in dimension five and in higher dimensions.

\begin{theorem}[Dimension $N=5$]\label{thm:main1}
Let $N=5$ and $k\ge2$.
\begin{enumerate}[label=\textup{(\roman*)}]
\item There is $\e_0>0$ such that for every $\e\in(0,\e_0)$ problem \eqref{pb} in $\mathcal B_1$ has a $\Gamma_k$-invariant solution $u_\e$ of the form \eqref{eq:form}, with $\|\phi_\e\|_{H^1_0}=O(\e^{5/2})$.
\item There is a constant $c_0>0$, independent of $k$ and $\e$, such that the radius $r_\e$ of the polygon satisfies $k\,(1-r_\e^2)<c_0$ for $\e$ small. In particular, every limit point $r_k$ of $r_\e$ as $\e\to0$ satisfies
\[
1-r_k^2<\frac{c_0}{k},\qquad\text{so that}\qquad r_k\to1\quad\text{as }k\to+\infty .
\]
\end{enumerate}
\end{theorem}
\begin{remark}
Part (ii) says that, in dimension five, the ring of negative bubbles is pushed towards the boundary as the number of peaks increases: the distance from $\partial\mathcal B_1$ is at most comparable with the mutual distance $\sim 2\pi/k$ of two consecutive negative peaks. Numerically the function $g$ has exactly one zero $r_k$ for every $k\ge2$, it is simple, and $k(1-r_k)$ seems to converge to a positive constant close to $4$, which suggests that the rate $1/k$ in (ii) is optimal; see Remark~\ref{rem:N5num}.\end{remark}

\begin{theorem}[Dimension $N\ge 6$]\label{thm:main2}
Let $N\ge6$ and $k\ge2$.
\begin{enumerate}[label=\textup{(\roman*)}]
\item If $2\le k\le6$, then for $\e$ small problem \eqref{pb} in $\mathcal B_1$ has two distinct $\Gamma_k$-invariant solutions $u^\pm_\e$ of the form \eqref{eq:form}, with radii $r^-_\e<2^{-1/2}<r^+_\e$.
\item There is $K_N\ge7$, depending only on $N$, such that for every $k\ge K_N$ the following holds: for every compact set $\mathcal K\subset\mathcal M$ there is $\e_1>0$ such that, for $\e\in(0,\e_1)$, problem \eqref{pb} has no solution of the form \eqref{eq:form} with $(r,d_0,d_1)\in\mathcal K$ and with $\phi_\e$ satisfying the orthogonality conditions and the estimate provided by the Lyapunov--Schmidt reduction \textup{(}see Corollary~\ref{cor:exist}\,\textup{(c)} for the precise statement\textup{)}.
\end{enumerate}
\end{theorem}

\begin{remark}
The restriction $k\le 6$ in (i) comes from the explicit computation of $g(2^{-1/2})$ (Proposition~\ref{prop:k6}); numerically, the sign change of $g$ persists for some larger values of $k$ in low dimensions (e.g.\ up to $k=17$ for $N=6$, $k=10$ for $N=7$, $k=8$ for $N=8$, $k=7$ for $N=10$ and so on). For $13\le N\le 30$, instead, numerically it occurs exactly for $2\le k\le 6$, while for $N\ge 31$ the value $k=7$ seems to appear again (the sign change takes place near $r\approx0.70$, slightly beyond $2^{-1/2}$, where $g$ is still positive).
\end{remark}
\medskip

\begin{figure}[htbp]
    \centering
    \includegraphics[width=0.7\textwidth]{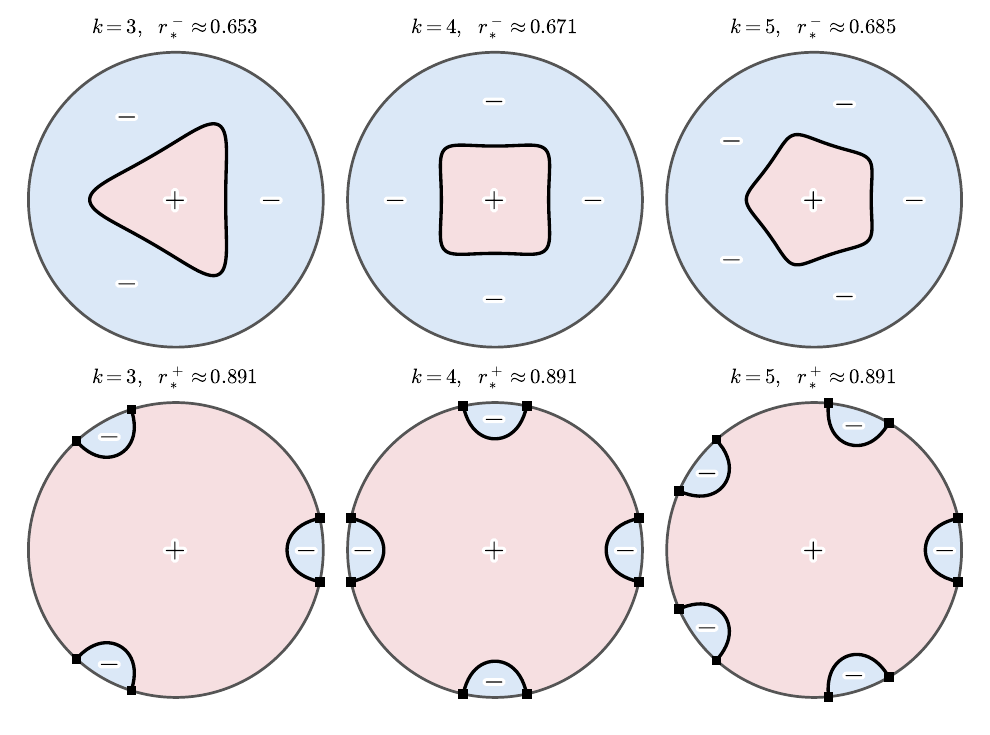}
    \caption{\tiny{Section in the plane $\{x_3=\dots=x_N=0\}$ containing the peaks of the limit profile $G(\cdot,0)-\sigma(r_*)\sum_jG(\cdot,\xi_j(r_*))$ of $\mu^{-\frac{N-2}{2}}u_\e$, for the two solutions $u^-_\e$ (top) and $u^+_\e$ (bottom) of Theorem~\ref{thm:main2} with $k=3,4,5$ peaks in dimension $N=6$. The positive region is shaded in red, the negative one in blue, and the squares mark the points where the nodal set (black curve) meets $\partial\mathcal B_1$. In this plane the nodal set of $u^-_\e$ is a closed curve around the origin which does not touch the boundary, while that of $u^+_\e$ consists of $k$ arcs which meet the boundary near the negative peaks. This does not contradict Theorem~\ref{thm:nodal}: outside this plane, near the points $(0,0,\pm1,0,\dots,0)$ of the sphere $\partial\mathcal B_1\cap\{x_1=x_2=0\}$, $u^-_\e$ is positive, so that its nodal set reaches $\partial\mathcal B_1$ there (see Remark~\ref{rem:two}).}}
    \label{fig:nodal-N6}
\end{figure}

We then study the nodal properties of these solutions: for $2\le k\le 6$ the exterior normal derivative changes sign on $\partial\mathcal B_1$ for \emph{every} solution of the form \eqref{eq:form}, so that the nodal set always reaches the boundary; for $N\ge6$ the two solutions of Theorem~\ref{thm:main2}\,(i) do so in very different ways (see Remark~\ref{rem:two} and Figure~\ref{fig:nodal-N6}).

\begin{theorem}\label{thm:nodal}
Let $N\ge5$ and $2\le k\le6$, and let $u_\e$ be a family of solutions of \eqref{pb} of the form \eqref{eq:form} given by Theorem~\ref{thm:main1} or by Theorem~\ref{thm:main2}\,\textup{(i)}. Then there is $\e_2>0$ such that for $\e\in(0,\e_2)$:
\begin{enumerate}[label=\textup{(\alph*)}]
\item $\partial_\nu u_\e<0$ on $\partial \mathcal B_1\cap\{x_1=x_2=0\}$; hence $u_\e>0$ in $\mathcal B_1$ near this $(N-3)$-dimensional sphere;
\item $\partial_\nu u_\e>0$ at the points $\xi_j(r_\e)/r_\e$, $j=1,\dots,k$ (which do not depend on $\e$); hence $u_\e<0$ in $\mathcal B_1$ near these points;
\item the nodal set $\mathcal N(u_\e):=\overline{\{x\in \mathcal B_1:\ u_\e(x)=0\}}$ intersects $\partial \mathcal B_1$.
\end{enumerate}
\end{theorem}

Moreover, by an energy argument, for every $k\ge2$ the solutions of Theorems~\ref{thm:main1} and~\ref{thm:main2}\,(i) have exactly one positive nodal domain and at most $k$ negative ones (Proposition~\ref{prop:domains}).

\begin{figure}[htbp]
    \centering
    \includegraphics[width=0.7\textwidth]{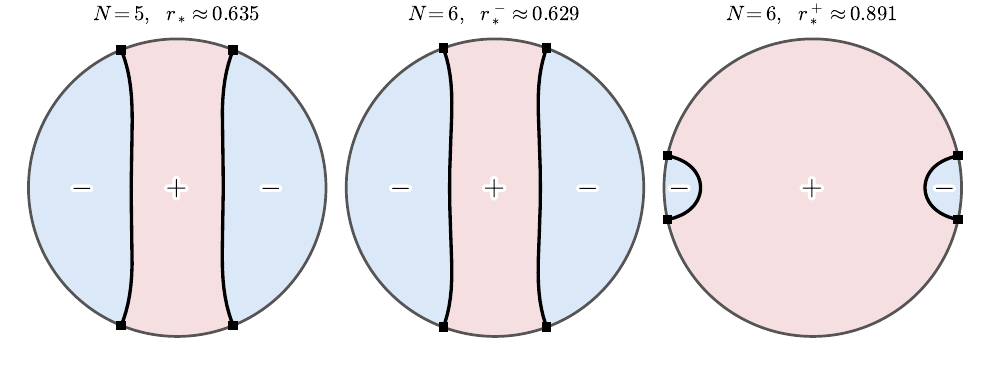}
    \caption{\tiny{The case $k=2$: the nodal set (black curve) in the plane $\{x_3=\dots=x_N=0\}$ of the limit profile $G(\cdot,0)-\sigma(r_*)\big(G(\cdot,\xi_1(r_*))+G(\cdot,\xi_2(r_*))\big)$ of $\mu^{-\frac{N-2}{2}}u_\e$, for $N=5$ (left) and for the two solutions $u^-_\e$, $u^+_\e$ with $N=6$ (centre and right). The positive region is shaded in red, the negative one in blue; the squares mark the points where the nodal set meets $\partial\mathcal B_1$, i.e.\ $x_1=\pm c_*$ (see Remark~\ref{rem:BDP}). Since for $k=2$ the limit profile depends only on $x_1$ and $|(x_2,\dots,x_N)|$, this section describes the whole nodal set.}}
    \label{fig:k2}
\end{figure}

\begin{remark}[Comparison with \cite{BDP}]\label{rem:BDP}
In \cite{BDP} Bartsch, D'Aprile and Pistoia studied the profile of sign-changing solutions of an almost critical problem in the ball, and showed (see \cite[Theorem~1.1\,(d)]{BDP}) that, for the configuration of three aligned bubbles with alternating signs (a positive one at the centre and two negative ones symmetric with respect to it), the reduced problem has two critical points and one of the two corresponding solutions has a nodal surface which does not touch the boundary. Our configuration with $k=2$ is the Brezis--Nirenberg counterpart of that one, and Theorem~\ref{thm:nodal} shows that the outcome is different. In both problems the sign of $\partial_\nu u_\e$ on $\partial\mathcal B_1$ is determined, as $\e\to0$, by the sign of $\sigma(r_*)P_{r_*}-1$, where $P_{r}$ is the boundary density \eqref{eq:P}, $\sigma$ is the limit ratio of the scales defined in \eqref{eq:s} and $r_*$ is the limit radius; a nodal surface not touching $\partial\mathcal B_1$ requires $\min_{\partial\mathcal B_1}\sigma(r_*)P_{r_*}\ge1$. In the Brezis--Nirenberg problem, however, the limit radius is a zero of the function $g$, and at \emph{every} zero of $g$ we have $\min\sigma P<1<\max\sigma P$ (Proposition~\ref{prop:mM}). Hence the nodal set of each of our solutions, including both solutions $u^\pm_\e$ of Theorem~\ref{thm:main2}\,(i) when $N\ge6$, reaches $\partial\mathcal B_1$: in particular, it cannot consist only of a closed surface around the positive bubble, as for one of the solutions in \cite{BDP}. The difference comes from the coupling between the scales and the location of the bubbles, which in our problem is governed by $g$ and forces the limit radius into the region where $\min\sigma P<1$.

For $k=2$ the boundary density depends only on $x_1$ and is increasing in $|x_1|$, so the limit nodal set meets $\partial\mathcal B_1$ along the two parallel $(N-2)$-dimensional spheres $\{x_1=\pm c_*\}$, where $c_*\in(0,1)$ solves $\sigma(r_*)P_{r_*}=1$. Numerically, $c_*\approx0.38$ for $N=5$, while for $N=6$ we find $c_*\approx0.35$ for $u^-_\e$ and $c_*\approx0.98$ for $u^+_\e$: in the latter case the negative bubbles are much more concentrated than the positive one and the negative nodal domains are thin regions near the boundary points $(\pm1,0,\dots,0)$; see Figure~\ref{fig:k2}.
\end{remark}

\begin{remark}[The two solutions for $N\ge6$]\label{rem:two}
The nodal sets of both solutions $u^\pm_\e$ of Theorem~\ref{thm:main2}\,(i) reach $\partial\mathcal B_1$, but for $3\le k\le6$ they do so in two very different ways. Let $\mathcal X:=\{x_3=\dots=x_N=0\}$ be the plane containing the peaks, and recall that $r^-_\e<2^{-1/2}<r^+_\e$. The following description is based on the analysis of the limit profile carried out in Remark~\ref{rem:twosol}, which is partly numerical.

For $u^-_\e$, the negative part of $u^-_\e$ near $\partial\mathcal B_1$ is a band around the whole circle $\partial\mathcal B_1\cap\mathcal X$, while the positive part is a neighbourhood of the $(N-3)$-dimensional sphere $\partial\mathcal B_1\cap\{x_1=x_2=0\}$; the nodal set meets $\partial\mathcal B_1$ along a set homeomorphic to $\mathbb S^1\times\mathbb S^{N-3}$ which separates them. In particular, \emph{in the plane $\mathcal X$ the nodal set of $u^-_\e$ is a closed curve around the origin which does not touch $\partial\mathcal B_1$}: the nodal set reaches the boundary only outside the plane of the peaks. Moreover, the $k$ negative peaks and the band around $\partial\mathcal B_1\cap\mathcal X$ lie in the same negative nodal domain, so that $u^-_\e$ has exactly two nodal domains.

For $u^+_\e$ the picture is the opposite one, at least when $k\in\{3,4\}$, or $k=5$ and $N\le9$, or $k=6$ and $N\le7$. The negative peaks are much more concentrated than the positive one (for $N=6$ the ratio of the scales is $s(r^+)\approx0.03$, against $s(r^-)\in(0.47,0.87)$), the negative part of $u^+_\e$ near $\partial\mathcal B_1$ consists of $k$ small disjoint caps centred at the points $\xi_j/|\xi_j|$, and the nodal set meets $\partial\mathcal B_1$ along $k$ disjoint $(N-2)$-dimensional spheres, one around each of these points. In the plane $\mathcal X$ the nodal set consists of $k$ arcs, each of which surrounds one negative peak and meets $\partial\mathcal B_1$ at two points, and $u^+_\e$ has $k$ negative nodal domains, i.e.\ the maximal number $k+1$ of nodal domains allowed by Proposition~\ref{prop:domains}. For larger values of $k$ and $N$, the nodal structure of $u^+_\e$ is the same as that of $u^-_\e$.

Thus, in the plane of the peaks, the nodal set of $u^-_\e$ looks like the closed nodal surface found in \cite{BDP}, although in $\mathcal B_1$ it does reach the boundary, outside the plane of the peaks; see Figure~\ref{fig:nodal-N6}. For $k=2$, instead, the two solutions have the same nodal structure, with three nodal domains, and differ only in the position of the two spheres where the nodal set meets $\partial\mathcal B_1$ (Remark~\ref{rem:BDP}). The same dichotomy holds for the solution of Theorem~\ref{thm:main1}: when $N=5$ and $3\le k\le6$ it has the structure of $u^-_\e$.
\end{remark}


\begin{remark}
In a forthcoming paper \cite{RAGO} other configurations that were open problems in \cite{PRV2} are considered: indeed,  it is shown that two other multi-peak configurations found in four dimensions in \cite{PRV2} survive in higher dimensions. In an annulus, for every even $k \geq 2$, solutions which blow up at the vertices of a regular $k$-gon with alternating signs exist. Instead, in the unit ball, solutions which blow up positively at the vertices of a regular $k$-gon lying in a coordinate plane and negatively at the vertices of a regular $k$-gon lying in an orthogonal plane, for $2\leq k\leq5$ if $N = 5, 6$ and for $2\leq k \leq4$ if $N \geq 7$, are constructed. Moreover, it is also shown that this range is sharp: for larger $k$ no solution of this type exists with parameters in a compact set, while, for $k$ in the above range the reduced energy has a unique critical point, which is non-degenerate.
\end{remark}

\subsection{Open problems and further directions}
The analysis above leaves several questions open, some of which are suggested by the numerical computations.

\begin{enumerate}[label=\textup{(\arabic*)}]
\item \emph{Uniqueness in dimension five.} Numerically, for $N=5$ the function $g$ has exactly one zero $r_k$ for every $k\ge2$, and this zero is simple (Remark~\ref{rem:N5num}). A proof of this fact would show that the solution of Theorem~\ref{thm:main1} comes from a non-degenerate critical point of Morse index one of the reduced function (Remark~\ref{rem:morse}) and would give its local uniqueness. It would also be interesting to prove that the rate $1/k$ in Theorem~\ref{thm:main1}\,(ii) is optimal, i.e.\ that $k(1-r_k)$ converges to a positive constant, as the numerics suggest.
\item \emph{The sharp range of $k$ when $N\ge6$.} Theorem~\ref{thm:main2} gives two solutions for $2\le k\le6$ and no solution for $k\ge K_N$, but the range $7\le k<K_N$ remains open. The numerical computations show that the set of $k$ for which $g$ changes sign depends on $N$ in a non-trivial way; it would be interesting to determine it exactly, and to prove that, whenever $g$ changes sign, it has exactly two zeros, both simple.
\item \emph{The nodal structure.} Theorem~\ref{thm:nodal} is proved for $2\le k\le6$; it is natural to ask whether the nodal set reaches $\partial\mathcal B_1$ for every $k$, in particular for $N=5$ and $k\ge7$. Moreover, by Proposition~\ref{prop:domains} our solutions have exactly one positive nodal domain and at most $k$ negative ones. Numerically, for $3\le k\le6$ the solution $u^-_\e$ has only one negative nodal domain, while $u^+_\e$ has $k$ of them, at least for $k\in\{3,4\}$ (see Remark~\ref{rem:two}); a rigorous proof of this dichotomy would be interesting.
\item \emph{Other domains and dimensions.} Our approach relies on the symmetries of the ball, which reduce the problem to three parameters. We believe, however, that this kind of solutions exists in every smooth bounded domain, at least in dimension five: this is the content of Conjecture~\ref{conj} below.
\end{enumerate}

\begin{conjecture}\label{conj}
Let $\Omega\subset\R^5$ be a smooth bounded domain and let $\xi_0\in\Omega$ be a minimum point of the Robin function $\tau$. Then there is $k_0=k_0(\Omega)\ge2$ such that, for every $k\ge k_0$ and every $\e$ small enough, problem \eqref{pb} has a sign-changing solution of the form
\[
u_\e=PU_{\delta_0^\e,\xi_0^\e}-\sum_{j=1}^kPU_{\delta_j^\e,\xi_j^\e}+\phi_\e,\qquad \delta_j^\e=d_j^\e\,\e,\qquad \|\phi_\e\|_{H^1_0(\Omega)}\to0,
\]
where, as $\e\to0$ and up to a subsequence, $\xi_0^\e\to\xi_0^{(k)}$ and $\xi_j^\e\to\xi_j^{(k)}$, $j=1,\dots,k$, with $\xi_0^{(k)}$, $\xi_1^{(k)},\dots,\xi_k^{(k)}$ distinct points of $\Omega$ such that
\[
\xi_0^{(k)}\to\xi_0\qquad\text{and}\qquad\max_{1\le j\le k}\mathrm{dist}\big(\xi_j^{(k)},\partial\Omega\big)=O\Big(\frac1k\Big)\qquad\text{as }k\to+\infty .
\]
In other words, the positive peak concentrates near the minimum point of the Robin function, while the negative peaks lie close to the boundary.
\end{conjecture}

Conjecture~\ref{conj} is suggested by the mechanism behind Theorem~\ref{thm:main1}. In the ball, by Proposition~\ref{prop:N5ring}, the negative peaks lie at distance $O(1/k)$ from the boundary, where their mutual interaction balances their interaction with $\partial\mathcal B_1$; moreover, they are much more concentrated than the positive peak. Indeed, the bounds $\sigma\le\sigma_0$ and \eqref{eq:betaa} established in the proof of Lemma~\ref{lem:layer} hold for every $N\ge5$, and for $N=5$ they give $k\beta(r_k)\sigma(r_k)=O(k^{-4})$ at every zero $r_k$ of $g$. Hence, by \eqref{eq:d0}, the scale of the positive peak converges, as $k\to+\infty$, to that of a single bubble, and its interaction with the ring becomes negligible: at main order the positive peak behaves as an isolated bubble, which concentrates at a critical point of the Robin function \cite{R}. In a general domain we expect the same balance to place the negative peaks in a thin layer near $\partial\Omega$, while the positive peak, almost decoupled from them, should sit close to a minimum point of $\tau$, as the centre does in the ball. For $N\ge6$ the situation is different, since by Theorem~\ref{thm:main2}\,\textup{(ii)} in the ball this configuration does not exist when $k$ is large.

\subsection{Organization of the paper}
The paper is organized as follows. In Section~\ref{sec:prelim} we recall the Lyapunov--Schmidt reduction of \cite{MiPi,MP,BMP} and reduce the problem to the study of the critical points of the function $\Psi$ of three variables defined in \eqref{eq:Psi}. In Section~\ref{sec:Psi} we prove the factorisation of $\nabla\Psi$ (Lemma~\ref{lem:factor}), we characterise the critical points of $\Psi$ and their local degree in terms of the function $g$, and we study the behaviour of $g$ at the endpoints of $(0,1)$. Sections~\ref{sec:N5} and~\ref{sec:k6} contain the proofs of Theorem~\ref{thm:main1}\,(i) and of Theorem~\ref{thm:main2}\,(i); the latter relies on an explicit formula for $g(2^{-1/2})$. In Section~\ref{sec:lower} we establish lower bounds for $g$, which give the non-existence result of Theorem~\ref{thm:main2}\,(ii) for large $k$ and the localisation of the ring of negative bubbles stated in Theorem~\ref{thm:main1}\,(ii). Finally, in Section~\ref{sec:nodal} we study the normal derivative of the solutions on $\partial\mathcal B_1$, we prove Theorem~\ref{thm:nodal} and we estimate the number of nodal domains.

\section{Preliminaries}\label{sec:prelim}

First we recall the definition of stable critical point used in this paper.
\begin{definition} \label{def:stable}
Let $f:D \to \R$ be a $C^1$-function, where $D \subset \R^n$ is an open set. We say that
$\xi_0$ is a stable critical point of $f$ if $\nabla f(\xi_0)=0$ and there exists a neighbourhood $V\subset\subset D$ of $\xi_0$ such that $\nabla f(\xi)\neq 0$ for any $\xi\in\partial V$,
\[\nabla f(\xi)=0,\quad \xi\in V \Longrightarrow f(\xi)=f(\xi_0)\]
and
\[\deg(\nabla f, V , 0)\neq  0,\] where $\deg$ denotes the Brouwer degree.
\end{definition}

Next we define the \emph{bubbles} (see \cite{A, CGS, T})
\begin{equation}\label{eq:bubble}
U_{\delta, \xi}(x)=\alpha_N \frac{\delta^{\frac{N-2}{2}}}{\left(\delta^2+|x-\xi|^2\right)^{\frac{N-2}{2}}},\quad x, \xi\in\R^N,\ \delta>0,\end{equation}
where $\alpha_N=\left[N(N-2)\right]^{\frac{N-2}{4}}$. As is well known, they are all the positive solutions of
\begin{equation}\label{eq:critico}-\Delta U=U^{\frac{N+2}{N-2}}\quad\text{in}\ \R^N.\end{equation}
We denote by $PU_{\delta, \xi}$ the projection of the bubble on $H^1_0(\Omega)$, i.e.
\[\Delta PU_{\delta, \xi}=\Delta U_{\delta, \xi}\quad \text{in}\ \Omega,\qquad PU_{\delta, \xi}=0\quad\text{on}\ \partial\Omega,\]
and we set $U:=U_{1,0}$ and
\begin{equation}\label{eq:aNbN}
a_N:=\int_{\R^N}U^{2^*-1}(x)\, dx,\qquad b_N:=\frac 12 \int_{\R^N}U^2(x)\, dx,\end{equation}
which are finite since $N\ge 5$.
In \cite{MiPi} the following result is proved.

\begin{theorem}[{\cite[Theorem 3.4]{MiPi}}]\label{thm:MiPi}
Let $(\bd d_0,\bd\xi_0)\in(0, +\infty)^k\times\Omega_k^*$ be a stable critical point of
$F^{\bd\sigma}_k$. Then, for $\e$ small enough, problem \eqref{pb} has a solution
\[
u_\e=\sum_{j=1}^k\sigma_jPU_{\delta_j^\e,\xi_j^\e}+\phi_\e
\]
such that $\delta_j^\e=d_j^\e\,\e^{\frac1{N-4}}$ and $\|\phi_\e\|_{H^1_0(\Omega)}\to0$ as $\e \to0$;
moreover, up to a subsequence, $d_j^\e\to d_j^*$ and $\xi_j^\e\to\xi_j^*$, where $(\bd d^*,\bd\xi^*)$ is a
critical point of $F^{\bd\sigma}_k$ with $F^{\bd\sigma}_k(\bd d^*,\bd\xi^*)=F^{\bd\sigma}_k(\bd d_0,\bd\xi_0)$.
\end{theorem}
Theorem~\ref{thm:MiPi} provides a general criterion for the existence of sign-changing solutions of
\eqref{pb} in dimension $N\ge5$; what remains is to find domains and sign patterns for which
$F^{\bd\sigma}_k$ has a stable critical point. In \cite{MiPi} this is done in two
situations. The first one is a dumbbell-shaped domain with thin handles, in which
problem \eqref{pb} has many sign-changing solutions with many nodal domains.
\begin{theorem}[{\cite[Theorem 1.4]{MiPi}}]\label{thm:MiPi1}
Let $N\ge 5$ and let $k\ge 2$ be a fixed integer. There exists a contractible domain $\Omega$
such that, for $\e$ small enough, problem \eqref{pb} has at least $2^{k-1}-1$
pairs of sign-changing solutions which blow up positively at $h$ points and negatively at $k-h$ points, for any
$h= 1, \ldots, k-1$.\end{theorem}
The second one is the ball, with peaks of alternating signs.
\begin{theorem}[{\cite[Theorem 1.5]{MiPi}}]\label{thm:MiPi2}
Let $\Omega$ be the unit ball in $\R^N$, $N\ge 5$, and let $h\ge1$ be a fixed integer. Then, for $\e$ small enough,
problem \eqref{pb} has at least a pair of sign-changing solutions which blow up positively at $h$ points
and negatively at $h$ points.\end{theorem}
In both cases the symmetries are used to reduce the problem to a function in which the location and the scale variables are not fully coupled.

From now on $\Omega=\mathcal B_1$ is the unit ball of $\R^N$. Solutions of \eqref{pb} are the critical points of
\[
J_\e(u)=\frac12\int_{\mathcal B_1}|\nabla u|^2-\frac\e2\int_{\mathcal B_1}u^2-\frac1{2^*}\int_{\mathcal B_1}|u|^{2^*},\qquad u\in H^1_0(\mathcal B_1).
\]
In the ball the Robin function and the Green function are explicitly given by
\begin{equation}\label{eq:GR}
\tau(x)=\frac{1}{(1-|x|^2)^{N-2}},\qquad
G(x, y)=\frac{1}{|x-y|^{N-2}}-\frac{1}{(1-2 \langle x,y \rangle + |x|^{2}|y|^{2})^{\frac{N-2}{2}}},
\end{equation}
up to the multiplicative constant $\omega_N^{-1}$, which we drop in order to simplify the computations.

\subsection{The configuration} Let $\xi_j(r)$, $r\in(0,1)$, be the points defined in
\eqref{eq:xi}. Since $J_\e$ is $\Gamma_k$-invariant, by the principle of symmetric criticality
the critical points of $J_\e|_{\mathcal H_k}$ are critical points of $J_\e$. We look for
solutions of the form
\begin{equation}\label{eq:form1}
u_\e=V_y+\phi,\qquad V_y:=PU_{\delta_0,0}-\sum_{j=1}^kPU_{\delta_1,\xi_j(r)},\qquad \delta_i=d_i\,\e^{\frac1{N-4}},\quad i=0,1,
\end{equation}
where $y=(r,d_0,d_1)$ lies in a compact subset of $\mathcal M$, and $\phi\in\mathcal H_k$ satisfies suitable orthogonality conditions and is such that
$\|\phi\|_{H^1_0}\to 0$ as $\e\to0$. For $\eta\in(0,\frac12)$ we set
\[\mathcal P_\eta:=[\eta, 1-\eta]\times[\eta, \eta^{-1}]^2\subset\mathcal M .\]
The Lyapunov--Schmidt reduction is by now standard and we do not repeat it here: we refer to
\cite{MiPi, MP, BMP} for the details and only recall its outcome. For every $y\in\mathcal P_\eta$ the reduction provides first the remainder $\phi_y\in\mathcal H_k$.
Because of the $\Gamma_k$-symmetry, the relevant kernel is three-dimensional:
\[
\Zs:=\mathrm{span}\Big\{\partial_\delta PU_{\delta,0}\big|_{\delta=\delta_0},\ \sum_{j=1}^k\partial_\delta PU_{\delta,\xi_j(r)}\big|_{\delta=\delta_1},\ \sum_{j=1}^k\partial_rPU_{\delta_1,\xi_j(r)}\Big\}\subset\mathcal H_k ;
\]
we denote by $\Zs^\perp$ its orthogonal complement in $\mathcal H_k$ and by $\Pi^\perp$ the orthogonal projection on it. Finally, let $\mu:=\e^{\frac{1}{N-4}}$ and
\[
\vartheta(\mu):=
\begin{cases}
\mu^{5/2}, & N=5,\\[1mm]
\mu^{4}|\log\mu|^{2/3}, & N=6,\\[1mm]
\mu^{\frac{N+2}{2}}, & N\ge 7.
\end{cases}
\]

\begin{proposition}[{\cite{MP,BMP}}]\label{prop:phi}
There are $\e_0,C>0$ such that for $\e\in(0,\e_0)$ and every $y\in\mathcal P_\eta$ there is a unique $\phi_y\in\Zs^\perp$ with $\|\phi_y\|\le C\vartheta(\mu)$ such that
\[
\Pi^\perp\Big[V_y+\phi_y-(-\Delta)^{-1}\big(|V_y+\phi_y|^{2^*-2}(V_y+\phi_y)+\e(V_y+\phi_y)\big)\Big]=0 .
\]
Moreover, the map $y\mapsto\phi_y$ is of class $C^1$.
\end{proposition}

Hence the whole problem is reduced to the study of a finite-dimensional energy, described by the following results.

\begin{lemma}[{\cite{MP,BMP}}]\label{lem:red}
Let $\widetilde J_\e(y):=J_\e(V_y+\phi_y)$. If $y\in\mathrm{int}\,\mathcal P_\eta$ is a critical point of $\widetilde J_\e$, then $V_y+\phi_y$ is a solution of \eqref{pb}.
\end{lemma}

\begin{proposition}[{\cite{R,MP,BMP}}]\label{prop:exp}
There are constants $A,B>0$, depending only on $N$, such that, as $\e\to0$,
\[
\widetilde J_\e(y)=(k+1)\frac{S^{N/2}}N+\mu^{N-2}\big(\Psi_{A,B}(y)+\varpi_\e(y)\big),\qquad \varpi_\e\to0\ \text{in }C^1(\mathcal P_\eta),
\]
where $S$ is the best Sobolev constant,
\begin{equation}\label{eq:PsiAB}
\Psi_{A,B}(r,d_0,d_1)=\frac A2\Big[d_0^{N-2}+k\varphi_k(r)d_1^{N-2}\Big]+kA\,\beta(r)\,(d_0d_1)^{\frac{N-2}{2}}-\frac B2\big(d_0^2+kd_1^2\big)
\end{equation}
and
\begin{equation}\label{eq:phibeta}
\beta(r):=G(0,\xi_1(r))=\frac{1}{r^{N-2}}-1,\qquad \varphi_k(r):=\tau(\xi_1(r))-\sum_{j=2}^kG(\xi_1(r),\xi_j(r)).
\end{equation}
\end{proposition}

Let $c:=(B/A)^{\frac{1}{N-4}}$. Then $\Psi_{A,B}(r,d_0,d_1)=Bc^2\,\Psi(r,d_0/c,d_1/c)$, where
\begin{equation}\label{eq:Psi}
\boxed{
\Psi(r,d_0,d_1)=\frac{1}{2}\Big[d_0^{N-2}+k\varphi_k(r)d_1^{N-2}\Big]+k\beta(r)(d_0d_1)^{\frac{N-2}{2}}-\frac{1}{2}\big(d_0^2+kd_1^2\big).
}
\end{equation}
Since this change of variables is a positive dilation in $(d_0,d_1)$, critical points and local degrees of $\nabla\Psi_{A,B}$ and $\nabla\Psi$ correspond. Hence, by Proposition~\ref{prop:exp} and Lemma~\ref{lem:red}, Theorems~\ref{thm:main1} and~\ref{thm:main2} reduce to the study of the critical points of $\Psi$, which is the content of the next section.

\section{The reduced function \texorpdfstring{$\Psi$}{Psi}}\label{sec:Psi}
\noindent Let $N \ge 5$ and $k\ge2$. Since $|\xi_1(r)-\xi_{j+1}(r)|=2r\sin\frac{\pi j}{k}$ and $\langle\xi_1(r),\xi_{j+1}(r)\rangle=r^2\cos\frac{2\pi j}{k}$, the explicit expressions \eqref{eq:GR} of the Green and Robin functions give
\[\begin{aligned}\varphi_k(r)
&=\frac{1}{(1-|\xi_1(r)|^2)^{N-2}}-\sum_{j=2}^k\left(\frac{1}{|\xi_1(r)-\xi_j(r)|^{N-2}}-\frac{1}{\left(1-2\langle\xi_1(r), \xi_j(r)\rangle+r^4\right)^{\frac{N-2}{2}}}\right)\\
&=\frac{1}{(1-r^2)^{N-2}}-\sum_{j=1}^{k-1}\left(\frac{1}{2^{N-2}r^{N-2}\sin^{N-2}\left(\frac{\pi j}{k}\right)}-\frac{1}{\left(1-2r^2\cos\left(\frac{2\pi j}{k}\right)+r^4\right)^{\frac{N-2}{2}}}\right).\end{aligned}\]
We set
 \[a=a(r):=1-r^2, \qquad \mathfrak s_j:=\sin\left(\frac{\pi j}{k}\right) \qquad\text{and}\qquad  b_j:=2r\mathfrak s_j,\qquad j=1,\dots,k-1.\]
Since
\[1-2r^2\cos\left(\frac{2\pi j}{k}\right)+r^4=a^2+b_j^2,\]
we can rewrite
\begin{equation}\label{eq:phiexpl}
\varphi_k(r)=a^{-(N-2)}-\sum_{j=1}^{k-1}\tau_j,\qquad \tau_j:=b_j^{-(N-2)}-(a^2+b_j^2)^{-\frac{N-2}{2}}>0.
\end{equation}
We remark that the functions $\varphi_k$ and $\beta$ are real analytic on $(0,1)$ and
\begin{equation}\label{eq:signs}
\beta(r)=r^{-(N-2)}-1>0,\qquad \beta'(r)=-(N-2)\,r^{-(N-1)}<0,\qquad\text{for all }r\in(0,1).
\end{equation}
In the following lemma we derive a useful expression of $\varphi_{k}'$ and we show that $\varphi_k$ is increasing on $(0, 1)$.
\begin{lemma}\label{lem:phiprime}
For every $r\in(0,1)$ we have
\begin{equation}\label{eq:phiprime}
\varphi_k'(r)=2(N-2)\,r\,a^{-(N-1)}+\sum_{j=1}^{k-1}\omega_j,
\end{equation}
where
\begin{equation}\label{eq:omegaj}
\omega_j=\frac {N-2}{r}\Big[b_j^{-(N-2)}-b_j^2(a^2+b_j^2)^{-\frac{N}{2}}\Big]+2(N-2)r\,a\,(a^2+b_j^2)^{-\frac{N}{2}}.
\end{equation}
In particular 
\[ \varphi_k'>0\quad\text{on}\quad (0,1).\]
\end{lemma}

\begin{proof}
Let us observe that 
\begin{equation*}
\begin{aligned}
\varphi'_{k}(r) & = 2r(N-2)a^{-(N-1)} - \sum_{j=1}^{k-1} \frac{d}{dr} \Big[b_{j}^{-(N-2)} - (a^{2}+b_{j}^{2})^{-\frac{N-2}{2}}\Big]
\\& = 2r(N-2)a^{-(N-1)} + \sum_{j=1}^{k-1} \Big[ (N-2) b_{j}^{-(N-1)} \frac{d b_j}{dr}- \frac{N-2}{2} (a^{2}+b_{j}^{2})^{-\frac{N}{2}} \frac{d}{dr}(a^{2}+b_{j}^{2})\Big]
\\&=2r(N-2)a^{-(N-1)} + \sum_{j=1}^{k-1} \Big[ \frac{2(N-2)\mathfrak s_{j}}{b_{j}^{N-1}} - \frac{N-2}{2}(a^{2}+b_{j}^{2})^{-\frac{N}{2}}(4b_{j}\mathfrak s_{j}-4ar)\Big],
\end{aligned}
\end{equation*}
where we used $\frac{d}{dr}b_j=2\mathfrak s_j=b_j/r$ and $\frac{d}{dr}(a^2+b_j^2)=4b_j\mathfrak s_j-4ar$. Hence \eqref{eq:phiprime} follows, with $\omega_j$ as in \eqref{eq:omegaj}.
Moreover, the bracket in \eqref{eq:omegaj} is positive, since $(a^2+b_j^2)^{-\frac{N}{2}}<b_j^{-N}$; therefore $\omega_{j}>0$ and
$\varphi_k'>0$ on $(0,1)$.
\end{proof}
\subsection{Critical points of $\Psi$} In order to look for critical points of $\Psi$, we compute the three partial derivatives of $\Psi$ with respect to its own variables $(r,d_0,d_1)$. Here we rewrite them in terms of the ratio $s=d_1/d_0$ and show that $\nabla\Psi$ factorises through three simpler functions $\mathcal R,\mathcal A,\mathcal G$. Then we solve $\mathcal R=\mathcal A=\mathcal G=0$ one equation at a time finding the critical point. At the end, we compute the sign of the Hessian determinant at a critical point and the local degree of $\nabla\Psi$.
\subsubsection{The partial derivatives of $\Psi$}\label{ssec:derivatives}

Differentiating term by term,
\begin{align}
\partial_r\Psi(r,d_0,d_1)&=\frac k2\,d_1^{N-2} \varphi_k'(r)+k\,(d_0d_1)^{\frac{N-2}{2}}\beta'(r),\label{eq:Psir}\\[1mm]
\partial_{d_0}\Psi(r,d_0,d_1)&=\frac{N-2}{2}\,d_0^{N-3}+\frac{N-2}{2}\,k\beta(r)\,d_1^{\frac{N-2}{2}}d_0^{\frac{N-4}{2}}-d_0,\label{eq:Psid0}\\[1mm]
\partial_{d_1}\Psi(r,d_0,d_1)&=\frac{N-2}{2}\,k\varphi_k(r)\,d_1^{N-3}+\frac{N-2}{2}\,k\beta(r)\,d_0^{\frac{N-2}{2}}d_1^{\frac{N-4}{2}}-k\,d_1.\label{eq:Psid1}
\end{align}
We now substitute $d_1=s\,d_0$, $s>0$, in \eqref{eq:Psir}--\eqref{eq:Psid1} obtaining
\begin{align}
\partial_r\Psi(r, d_0, d_1)&=k\,d_0^{N-2}s^{\frac{N-2}{2}}\Big[\tfrac12\varphi_k'(r)\,s^{\frac{N-2}{2}}+\beta'(r)\Big],\label{eq:Psir-s}\\[1mm]
\partial_{d_0}\Psi(r, d_0, d_1)&=d_0\Big[\tfrac{N-2}{2}\,d_0^{N-4}\big(1+k\beta(r)s^{\frac{N-2}{2}}\big)-1\Big],\label{eq:Psid0-s}\\[1mm]
\partial_{d_1}\Psi(r, d_0, d_1)&=k\,s\,d_0\Big[\tfrac{N-2}{2}\,d_0^{N-4}\left(\varphi_k(r)s^{N-4}+\beta(r)s^{\frac{N-6}{2}}\right)-1\Big].\label{eq:Psid1-s}
\end{align}
For $r\in(0,1)$, $d_0>0$ and $s>0$ we then introduce the functions
\begin{align}
\mathcal R(r,s)&:=\tfrac12\varphi_k'(r)\,s^{\frac{N-2}{2}}+\beta'(r),\label{eq:R}\\
\mathcal A(r,d_0,s)&:=\tfrac{N-2}{2}\,d_0^{N-4}\big(1+k\beta(r)s^{\frac{N-2}{2}}\big)-1,\label{eq:A}\\
\mathcal G(r,s)&:=1+k\beta(r)s^{\frac{N-2}{2}}-\varphi_k(r)s^{N-4}-\beta(r)s^{\frac{N-6}{2}}.\label{eq:G}
\end{align}
Note that $\mathcal R$ and $\mathcal G$ do not depend on $d_0$.

\begin{lemma}\label{lem:factor}
For every $(r,d_0,d_1)\in\mathcal M$, with $s=\frac{d_1}{d_0}$,
\begin{equation}\label{eq:factor}
\nabla\Psi(r,d_0,d_1)
=
\mathbb T(d_0,s)
\begin{pmatrix}
\mathcal R(r,s)\\
\mathcal A(r,d_0,s)\\
\mathcal G(r,s)
\end{pmatrix}
\end{equation}
with
\begin{equation}\label{eq:matrix-M}
\mathbb T(d_0,s)
=
\begin{pmatrix}
k\,d_0^{N-2}s^{\frac{N-2}{2}} & 0 & 0\\[1mm]
0 & d_0 & 0\\[1mm]
0 & k\,s\,d_0 & -\tfrac{N-2}{2}\,k\,s\,d_0^{N-3}
\end{pmatrix}
\end{equation}
and
\[
\det \mathbb T(d_0,s)=-\tfrac{N-2}{2}\,k^2\,s^{\frac{N}{2}}d_0^{2(N-2)}<0 .
\]
In particular $\mathbb T(d_0,s)$ is invertible on $\mathcal M$.
\end{lemma}
\begin{proof}
The first two rows are \eqref{eq:Psir-s} and \eqref{eq:Psid0-s}: 
\[\partial_r\Psi(r, d_0, d_1)=kd_0^{N-2} s^{\frac{N-2}{2}} \mathcal R\quad\text{and}\quad \partial_{d_0}\Psi(r, d_0, d_1)=d_0\mathcal A.\]
Indeed, by \eqref{eq:Psir} and $d_1=s\,d_0$,
\[
\partial_r\Psi=\frac k2\varphi_k's^{N-2}d_0^{N-2}+k\beta's^{\frac{N-2}2}d_0^{N-2}=kd_0^{N-2}s^{\frac{N-2}2}\Big[\tfrac12\varphi_k's^{\frac{N-2}2}+\beta'\Big],
\]
and by \eqref{eq:Psid0}
\[
\partial_{d_0}\Psi=\frac{N-2}2d_0^{N-3}+\frac{N-2}2k\beta s^{\frac{N-2}2}d_0^{N-3}-d_0=d_0\Big[\frac{N-2}2d_0^{N-4}\big(1+k\beta s^{\frac{N-2}2}\big)-1\Big].
\]
For the third row we compute
\[
\begin{aligned}
\mathcal A-\tfrac{N-2}{2}d_0^{N-4}\mathcal G
&=\tfrac {N-2}{2} d_0^{N-4}\Big[\big(1+k\beta s^{\frac{N-2}{2}}\big)-\big(1+k\beta s^{\frac{N-2}{2}}-\varphi_ks^{N-4}-\beta s^{\frac{N-6}{2}}\big)\Big]-1\\
&=\tfrac{N-2}{2}d_0^{N-4}\big(\varphi_ks^{N-4}+\beta s^{\frac{N-6}{2}}\big)-1,
\end{aligned}
\]
which is the bracket in \eqref{eq:Psid1-s}. Hence
\[
\partial_{d_1}\Psi(r, d_0, d_1)=ksd_0\Big(\mathcal A-\tfrac {N-2}{2}d_0^{N-4}\mathcal G\Big)=ksd_0\,\mathcal A-\tfrac {N-2}{2}ksd_0^{N-3}\,\mathcal G .
\]
Finally $\mathbb T$ is lower triangular, so 
\[\det \mathbb T=kd_0^{N-2}s^{\frac{N-2}{2}}\cdot d_0\cdot\left(-\frac{N-2}{2}ksd_0^{N-3}\right)=-\frac{N-2}{2}k^2s^{\frac{N}{2}}d_0^{2(N-2)}<0.\]
\end{proof}

\subsubsection{Solving $\nabla\Psi=0$}\label{ssec:solve}

By Lemma~\ref{lem:factor}, $\nabla\Psi(r,d_0,d_1)=0$ if and only if
\[
\mathcal R(r,s)=0,\qquad \mathcal A(r,d_0,s)=0,\qquad \mathcal G(r,s)=0,\qquad s=d_1/d_0 .
\]
We solve these three equations in this order.

\smallskip\noindent\emph{Step 1: the equation $\mathcal R=0$ determines $s$.}
By \eqref{eq:signs} we can define, for $r\in(0,1)$,
\begin{equation}\label{eq:s}
\sigma(r):=\frac{2|\beta'(r)|}{\varphi_k'(r)}=\frac{2(N-2)\,r^{-(N-1)}}{\varphi_k'(r)}>0\qquad \text{and}\qquad s(r):=\sigma(r)^{\frac{2}{N-2}}>0.
\end{equation}
Since $\beta'=-|\beta'|=-\frac12\varphi_k'\sigma$, we can rewrite \eqref{eq:R} as
\begin{equation}\label{eq:Rfactor}
\mathcal R(r,s)=\tfrac12\varphi_k'(r)\big(s^{\frac{N-2}{2}}-\sigma(r)\big)=\tfrac12\varphi_k'(r)\big(s^{\frac{N-2}{2}}-s(r)^{\frac{N-2}{2}}\big).
\end{equation}
As $\varphi_k'>0$ and $s\mapsto s^{\frac{N-2}{2}}$ is strictly increasing, $\mathcal R(r,s)$ has the sign of $s-s(r)$; in particular
\[\mathcal R(r,s)=0\iff s=s(r)\iff s^{\frac{N-2}{2}}=\sigma(r).\]
\smallskip\noindent\emph{Step 2: the equation $\mathcal G=0$ becomes an equation in $r$ only.}
Substituting $s=s(r)$ in \eqref{eq:G} we obtain the function
\begin{equation}\label{eq:g}
\begin{aligned}
g(r) :&=\mathcal G(r,s(r))=1+k\beta(r) s(r)^{\frac{N-2}{2}}-\varphi_k(r)s(r)^{N-4}-\beta(r) s(r)^{\frac{N-6}{2}}\\&=1+k\beta(r)\sigma(r)-s(r)^{-2}\sigma(r)\,(\varphi_k(r)\sigma(r)+\beta(r)).
\end{aligned}
\end{equation}

Thus, on $\{\mathcal R=0\}$,
\[
\mathcal G(r,s)=0\iff g(r)=0 .
\]

\smallskip\noindent\emph{Step 3: the equation $\mathcal A=0$ determines $d_0$.}
For $s=s(r)$, by \eqref{eq:signs} we have
\[1+k\beta s^{\frac{N-2}{2}}=1+k\beta\sigma>1,\]
and \eqref{eq:A} gives
\[\mathcal A(r,d_0,s(r))=\frac{N-2}{2} d_0^{N-4}(1+k\beta(r)\sigma(r))-1,\]
which is strictly increasing in $d_0$ (since $N \ge 5$) and vanishes exactly at
\begin{equation}\label{eq:d0}
d_0=d_0(r):=\left(\frac{\frac{2}{N-2}}{1+k\beta(r)\sigma(r)}\right)^{\frac{1}{N-4}}.
\end{equation}
The functions $\sigma$, $s$, $g$ and $d_0$ are real analytic on $(0,1)$, since $\varphi_k$ and $\beta$ are, $\varphi_k'>0$, $\sigma>0$ and $1+k\beta\sigma>0$. Summing up, we have proved the following result.
\begin{lemma}\label{lem:crit}
\begin{enumerate}[label=\textup{(\roman*)}]
\item $(r,d_0,d_1)\in\mathcal M$ is a critical point of $\Psi$ if and only if
\[
g(r)=0,\qquad d_0=d_0(r),\qquad d_1=s(r)\,d_0(r).
\]
\item At a critical point, $\sign(\det D^2\Psi(r,d_0,d_1))=\sign g'(r)$.
\end{enumerate}
\end{lemma}

\begin{proof}
(i) Since $\mathbb T$ is invertible (Lemma~\ref{lem:factor}), $\nabla\Psi=0$ if and only if $\mathcal R=\mathcal A=\mathcal G=0$. By Step~1, $\mathcal R=0$ if and only if $s=s(r)$, i.e.\ $d_1=s(r)d_0$;
then, by Step~2, $\mathcal G=0$ if and only if $g(r)=0$, and by Step~3, $\mathcal A=0$ if and only if $d_0=d_0(r)$.

(ii) Write $\mathcal F:=(\mathcal R,\mathcal A,\mathcal G)$, which is a function of $(r,d_0,s)$, and $\Theta(r,d_0,d_1):=(r,d_0,d_1/d_0)$, so that \eqref{eq:factor} reads
\[
\nabla\Psi=\mathbb T\,(\mathcal F\circ\Theta).
\]
Let $y_*=(r_*,d_{0*},d_{1*})$ be a critical point, so that $\mathcal F(\Theta(y_*))=0$. Differentiating the product, the terms containing derivatives of $\mathbb T$ are multiplied by $\mathcal F(\Theta(y_*))=0$; hence
\[
D^2\Psi(y_*)=D(\nabla\Psi)(y_*)=\mathbb T\cdot D_{(r,d_0,s)}\mathcal F(\Theta(y_*))\cdot D\Theta(y_*),
\]
where
\[
D\Theta=\begin{pmatrix}1&0&0\\0&1&0\\0&-d_1/d_0^2&1/d_0\end{pmatrix},\qquad\det D\Theta=\frac1{d_0},
\]
and therefore
\begin{equation}\label{eq:detHess}
\det D^2\Psi(y_*)=\det \mathbb T\cdot\det D_{(r,d_0,s)}\mathcal F\cdot\frac1{d_{0*}} .
\end{equation}
It remains to compute $\det D_{(r,d_0,s)}\mathcal F$. Since $\mathcal R$ and $\mathcal G$ do not depend on $d_0$,
\[
D_{(r,d_0,s)}\mathcal F=\begin{pmatrix}\mathcal R_r&0&\mathcal R_s\\ \mathcal A_r&\mathcal A_{d_0}&\mathcal A_s\\ \mathcal G_r&0&\mathcal G_s\end{pmatrix},\quad\text{and}\quad \det D_{(r,d_0,s)}\mathcal F=\mathcal A_{d_0}\,\big(\mathcal R_r\mathcal G_s-\mathcal R_s\mathcal G_r\big).
\]
From \eqref{eq:R}, we get
\[\mathcal R_s=\frac{N-2}{4}\varphi_k'(r)s^{\frac{N-4}{2}}>0.\]
From \eqref{eq:A}, we get
\[\mathcal A_{d_0}=\frac {N-2}{2} (N-4)\,d_0^{N-5}\big(1+k\beta(r) s^{\frac{N-2}{2}}\big)>0,\]
since $N\ge5$ and $\beta(r)>0$.
Finally, since $\mathcal R(r,s(r))\equiv0$ and $\mathcal R_s\ne0$, differentiating in $r$ gives
\[\mathcal R_r+\mathcal R_s\,s'(r)=0,\] 
i.e.\ $s'=-\frac{\mathcal R_r}{\mathcal R_s}$. Differentiating $g(r)=\mathcal G(r,s(r))$,
\[g'(r)=\mathcal G_r+\mathcal G_s\,s'(r)=\mathcal G_r-\mathcal G_s\frac{\mathcal R_r}{\mathcal R_s}
\quad\Longrightarrow\quad \mathcal R_r\mathcal G_s-\mathcal R_s\mathcal G_r=-\mathcal R_s\,g'(r),\]
all functions being evaluated at $(r,s(r))$.
Consequently, at $\Theta(y_*)$,
\[\det D_{(r,d_0,s)}\mathcal F=-\mathcal A_{d_0}\,\mathcal R_s\,g'(r_*),\qquad \mathcal A_{d_0}\mathcal R_s>0 .\]
Inserting the last identity and $\det \mathbb T<0$ (Lemma~\ref{lem:factor}) in \eqref{eq:detHess}, we easily get that
\[\det D^2\Psi(y_*)=\underbrace{\det \mathbb T}_{<0}\cdot\big(\underbrace{-\mathcal A_{d_0}\mathcal R_s}_{<0}\big)\,g'(r_*)\cdot\underbrace{\frac1{d_{0*}}}_{>0},\]
so $\sign(\det D^2\Psi(y_*))={\sign} g'(r_*)$.
\end{proof}

\subsubsection{A degree formula}\label{ssec:degree}
 
Lemma~\ref{lem:crit} describes the critical points of $\Psi$ and, at simple zeros of $g$,
their local degree. To handle zeros of $g$ that are possibly degenerate, we compute the
Brouwer degree of $\nabla\Psi$ on a whole region in terms of the sign of $g$ at two points.
 
\medskip\noindent\textbf{The set $\mathfrak U$.}
Fix $0<r_1<r_2<1$ with $g(r_1)\,g(r_2)\neq 0$. The functions $d_0(\cdot)$ and $s(\cdot)$
in \eqref{eq:d0}, \eqref{eq:s} are continuous and positive on $[r_1,r_2]$, so we can fix
$\eta\in(0,\frac12)$ such that
\begin{equation}\label{eq:eta}
  2\eta\le d_0(r)\le (2\eta)^{-1},\qquad 2\eta\le s(r)\le (2\eta)^{-1}
  \qquad\text{for all } r\in[r_1,r_2].
\end{equation}
In the coordinates $(r,d_0,s)$, with $s=d_1/d_0$, we work on the box
\[
\mathfrak  V:=(r_1,r_2)\times(\eta,\eta^{-1})\times(\eta,\eta^{-1}),
\]
and in the original coordinates on
\[
\mathfrak  U:=\Theta^{-1}(\mathfrak V)=\big\{(r,d_0,d_1)\,:\ r_1<r<r_2,\ \eta<d_0<\eta^{-1},\ \eta<d_1/d_0<\eta^{-1}\big\},
\]
where $\Theta(r,d_0,d_1)=(r,d_0,d_1/d_0)$ is a diffeomorphism of $\overline {\mathfrak U}$ onto $\overline {\mathfrak V}$.
Let
\[
  Z:=\big\{(r,d_0,s)\in\overline{\mathfrak V}\ :\ g(r)=0,\ d_0=d_0(r),\ s=s(r)\big\}.
\]
$Z$ is compact and, by $g(r_1)g(r_2)\ne0$ and \eqref{eq:eta}, it does not meet $\partial \mathfrak V$.
By Lemma~\ref{lem:crit}\,(i) the zeros of $\nabla\Psi$ in $\overline {\mathfrak U}$ are exactly $\Theta^{-1}(Z)$,
so $\nabla\Psi\neq0$ on $\partial \mathfrak U$ and $\deg(\nabla\Psi,\mathfrak U,0)$ is well defined.
 
\begin{proposition}\label{prop:degree}
Under the above assumptions,
\[
  \deg(\nabla\Psi,\mathfrak U,0)\;=\;\deg\big(g,(r_1,r_2),0\big)\;=\;\frac12\Big(\sign g(r_2)-\sign g(r_1)\Big).
\]
\end{proposition}
 
\begin{proof}
The second equality is the usual formula for the degree of a continuous real function on
an interval, so we only prove the first one. We deform $\nabla\Psi$, through admissible
homotopies, into the map
\[
 \mathcal  L(r,d_0,s):=\big(s-s(r),\ d_0-d_0(r),\ g(r)\big),
\]
whose degree can be computed by hand. Two elementary facts are used throughout.
 
\begin{itemize}
\item[(H)] If $H:[0,1]\times\overline {\mathfrak V}\to\R^3$ is continuous and $\{H(t,\cdot)=0\}=Z$ for every
$t\in[0,1]$, then $\deg(H(0,\cdot),\mathfrak V,0)=\deg(H(1,\cdot),\mathfrak V,0)$. Indeed $Z\cap\partial \mathfrak V=\emptyset$, so
the homotopy is admissible.
\item[(S)] If two real functions $f_1,f_2$ have the same sign at every point (with the same zeros),
then $(1-t)f_1+tf_2$ has that sign too and vanishes exactly where $f_1$ does, for every $t\in[0,1]$.
\end{itemize}
 
Recall that $\mathcal F=(\mathcal R,\mathcal A,\mathcal G)$ is a function of $(r,d_0,s)$ and that $\nabla\Psi=\mathbb T\,(\mathcal F\circ\Theta)$ by
Lemma~\ref{lem:factor}.
 
\smallskip\noindent\emph{Step 1: removing the matrix $\mathbb T$.} Here we show that \[\deg(\nabla\Psi,\mathfrak U,0)=\ -\deg(\mathcal F\circ\Theta,\mathfrak U,0)\]
Let $J:=\operatorname{diag}(1,1,-1)$ and $\mathbb T_t:=(1-t)\mathbb T+tJ$, $t\in[0,1]$. By \eqref{eq:matrix-M},
$\mathbb T_t(d_0,s)$ is lower triangular with diagonal entries
\[
  (1-t)\,k\,d_0^{N-2}s^{\frac{N-2}{2}}+t>0,\qquad (1-t)\,d_0+t>0,\qquad
  -(1-t)\tfrac{N-2}{2}\,k\,s\,d_0^{N-3}-t<0,
\]
hence $\det \mathbb T_t<0$ on $\overline{\mathfrak U}$ for every $t$. In particular $\mathbb T_t(\mathcal F\circ\Theta)$ vanishes exactly
on $\Theta^{-1}(Z)$, and by homotopy invariance
\[
  \deg(\nabla\Psi,\mathfrak U,0)=\deg\big(J\,(\mathcal F\circ\Theta),\mathfrak U,0\big)=\sign(\det J)\,\deg(\mathcal F\circ\Theta,\mathfrak U,0)=-\deg(\mathcal F\circ\Theta,\mathfrak U,0).
\]
 
\smallskip\noindent\emph{Step 2: change of variables.} Here we show that
\[-\deg(\mathcal F\circ\Theta,\mathfrak U,0)
=\ -\deg(\mathcal F,\mathfrak V,0)\]
Indeed, $\Theta:\overline {\mathfrak U}\to\overline{\mathfrak V}$ is a diffeomorphism with $\det D\Theta=1/d_0>0$, hence it preserves
orientation and $\deg(\mathcal F\circ\Theta,\mathfrak U,0)=\deg(\mathcal F,\mathfrak V,0)$.
 
\smallskip\noindent\emph{Step 3: from $ \mathcal F$ to $ \mathcal L$.} Here we show that \[\ -\deg(\mathcal F,\mathfrak V,0)
=\ -\deg(\mathcal L,\mathfrak V,0).\]

We modify one component at a time; by (H), at each stage it suffices to check that the zero set is $Z$.
 
\emph{(3a) First component.} By \eqref{eq:Rfactor} and $\varphi_k'>0$, $ \mathcal R(r,s)$ has the sign of
$s-s(r)$. By (S), along $\big((1-t) \mathcal R+t(s-s(r)),\, \mathcal A,\, \mathcal G\big)$ the first component vanishes iff
$s=s(r)$, so the zero set is always $Z$. Hence $ \mathcal F$ can be replaced by
\[
   \mathcal F_1:=\big(s-s(r),\  \mathcal A(r,d_0,s),\  \mathcal G(r,s)\big).
\]
 
\emph{(3b) Freezing $s=s(r)$ in $ \mathcal A$ and $ \mathcal G$.} For $t\in[0,1]$ let $s_t:=s(r)+t\,(s-s(r))>0$ and
\[
  H_t:=\big(s-s(r),\  \mathcal A(r,d_0,s_t),\  \mathcal G(r,s_t)\big),\qquad H_1=\mathcal F_1.
\]
A zero of $H_t$ has $s=s(r)$, hence $s_t=s(r)$, and there the last two components are
$ \mathcal A(r,d_0,s(r))$ and $ \mathcal G(r,s(r))=g(r)$, independently of $t$. So the zero set of $H_t$ is $Z$ for every $t$,
and $ \mathcal F_1$ can be replaced by
\[
  H_0=\big(s-s(r),\  \mathcal A(r,d_0,s(r)),\ g(r)\big).
\]
 
\emph{(3c) Second component.} By Step~3 of Subsection~\ref{ssec:solve}, $d_0\mapsto \mathcal A(r,d_0,s(r))$ is strictly
increasing and vanishes only at $d_0(r)$, i.e.\ it has the sign of $d_0-d_0(r)$. By (S) we can replace
it by $d_0-d_0(r)$, obtaining $\mathcal L$. Therefore $\deg(\mathcal F,\mathfrak V,0)=\deg(\mathcal L,\mathfrak V,0)$.
 
\smallskip\noindent\emph{Step 4: the degree of $\mathcal L$.} Here we show that \[
 -\deg(\mathcal L,\mathfrak V,0)
= \deg\big(g,(r_1,r_2),0\big).\] Indeed, the map $\Xi(r,d_0,s):=\big(r,\,d_0-d_0(r),\,s-s(r)\big)$ is a diffeomorphism of $\overline {\mathfrak V}$ onto
its image, with $\det D\Xi=1$, and
\[
  (\mathcal L\circ\Xi^{-1})(r,u,v)=\big(v,\ u,\ g(r)\big).
\]
Hence $\deg(\mathcal L,\mathfrak V,0)=\deg(\mathcal L\circ\Xi^{-1},\Xi(\mathfrak V),0)$. The zeros of $\mathcal L\circ\Xi^{-1}$ are the points $(r,0,0)$ with
$g(r)=0$; they lie in the box $W:=(r_1,r_2)\times(-\eta,\eta)^2$, which is contained in $\Xi(\mathfrak V)$ by
\eqref{eq:eta} (if $|u|<\eta$ then $\eta<d_0(r)+u<(2\eta)^{-1}+\eta<\eta^{-1}$, and the same for $v$).
By excision we can compute the degree on $W$. There
\[
  (v,u,g(r))=\mathbb P\,\big(g(r),u,v\big),\qquad \mathbb P(x_1,x_2,x_3):=(x_3,x_2,x_1),\quad \det \mathbb P=-1,
\]
and the map $(r,u,v)\mapsto(g(r),u,v)$ is a product, whose degree on $W$ is
$\deg\big(g,(r_1,r_2),0\big)\cdot\deg\big(\mathrm{id},(-\eta,\eta)^2,0\big)=\deg\big(g,(r_1,r_2),0\big)$. Therefore
\[
  \deg(\mathcal L,\mathfrak V,0)=-\deg\big(g,(r_1,r_2),0\big).
\]
Combining Steps 1--4 gives $\deg(\nabla\Psi,\mathfrak U,0)=\deg\big(g,(r_1,r_2),0\big)$.
\end{proof}
\begin{remark}
If $g$ has a single zero $r_*\in(r_1,r_2)$ and it is simple, Proposition~\ref{prop:degree} agrees
with Lemma~\ref{lem:crit}\,(ii): the local degree of the non-degenerate critical point is
\[
\sign\det D^2\Psi=\sign g'(r_*)=\tfrac12\big(\sign g(r_2)-\sign g(r_1)\big).
\]
\end{remark}
\begin{corollary}\label{cor:exist}
Let $N\ge5$ and $k\ge2$.
\begin{enumerate}[label=\textup{(\alph*)}]
\item If $g(r_1)g(r_2)<0$ for some $0<r_1<r_2<1$, then for $\e$ small problem \eqref{pb} has a $\Gamma_k$-invariant solution of the form \eqref{eq:form} with $r_\e\in(r_1, r_2)$; as $\e\to0$, $r_\e$ accumulates only at zeros of $g$ in $(r_1,r_2)$.
\item If $g>0$ near $0$ and near $1$, and $g(\hat r)<0$ for some $\hat r\in(0,1)$, then for $\e$ small there are two distinct such solutions, with $r^-_\e<\hat r<r^+_\e$.
\item If $g>0$ on $(0,1)$, then for every compact set $\mathcal K\subset\mathcal M$ there is $\e_1>0$ such that, for $\e\in(0,\e_1)$, problem \eqref{pb} has no solution of the form $V_y+\phi$ with $y\in\mathcal K$, $\phi\in\Zs^\perp$ and $\|\phi\|\le C\vartheta(\mu)$, where $C$ is the constant of Proposition~\ref{prop:phi}.
\end{enumerate}
\end{corollary}
\begin{proof}
Recall that
\[
  \Psi_{A,B}(r,d_0,d_1)=Bc^2\,\Psi(r,d_0/c,d_1/c)\] and that \[  \nabla\Psi_{A,B}(r,d_0,d_1)=Bc^2\operatorname{diag}(1,c^{-1},c^{-1})\,\nabla\Psi(r,d_0/c,d_1/c)
\]
In particular $y=(r,d_0,d_1)$ is a critical point of $\Psi_{A,B}$ iff $(r,d_0/c,d_1/c)$ is a critical point
of $\Psi$, and by Lemma~\ref{lem:crit}\,(i) this forces $g(r)=0$.
 
\smallskip\noindent\emph{(a)} Let $\mathfrak U$ and $\eta$ be as in Subsection~\ref{ssec:degree} and
$\mathfrak U_c:=\{(r,cd_0,cd_1):(r,d_0,d_1)\in \mathfrak U\}$. Since $\Psi_{A,B}$ is obtained from $\Psi$ by an
orientation-preserving dilation and multiplication by a positive diagonal matrix,
Proposition~\ref{prop:degree} gives
\[
  \deg(\nabla\Psi_{A,B},\mathfrak U_c,0)=\deg(\nabla\Psi,\mathfrak U,0)=\tfrac12\big(\sign g(r_2)-\sign g(r_1)\big)=\pm1,
\]
and $\nabla\Psi_{A,B}\ne0$ on $\partial \mathfrak U_c$.
 
On $\overline{\mathfrak U_c}$ we have $r\in[r_1,r_2]$, $d_0\in[c\eta,c\eta^{-1}]$ and $d_1=d_0\,(d_1/d_0)\in[c\eta^2,c\eta^{-2}]$.
We can therefore fix $\eta'\in(0,\frac12)$ with $\overline{\mathfrak U_c}\subset\interior\mathcal P_{\eta'}$, and apply
Propositions~\ref{prop:phi} and~\ref{prop:exp} with $\eta'$ in place of $\eta$. They give
\[
  \mu^{-(N-2)}\nabla\widetilde J_\e=\nabla\Psi_{A,B}+\nabla\varpi_\e
  \quad\text{on }\overline{\mathfrak U_c},\qquad \|\nabla\varpi_\e\|_{L^\infty(\overline{\mathfrak U_c})}\to0.
\]
Let $\mathfrak c:=\min_{\partial \mathfrak U_c}|\nabla\Psi_{A,B}|>0$. For $\e$ so small that $|\nabla\varpi_\e|<\mathfrak c$ on $\overline{\mathfrak U_c}$,
the homotopy $\nabla\Psi_{A,B}+t\nabla\varpi_\e$, $t\in[0,1]$, does not vanish on $\partial \mathfrak U_c$. Hence
\[
  \deg(\nabla\widetilde J_\e,\mathfrak U_c,0)=\deg(\nabla\Psi_{A,B},\mathfrak U_c,0)=\pm1,
\]
and $\widetilde J_\e$ has a critical point $y_\e=(r_\e,d_{0,\e},d_{1,\e})\in \mathfrak U_c$. By
Lemma~\ref{lem:red}, $u_\e:=V_{y_\e}+\phi_{y_\e}$ solves \eqref{pb}. It lies in $\mathcal H_k$, and
$\|\phi_{y_\e}\|\le C\vartheta(\mu)=o\big(\mu^{\frac{N-2}2}\big)$, so it is of the form \eqref{eq:form}.
 
Let $y_*$ be a limit point of $y_\e$ along $\e_n\to0$. Then $y_*\in\overline{\mathfrak U_c}$ and, by the
uniform convergence above, $\nabla\Psi_{A,B}(y_*)=0$. Thus $y_*\notin\partial \mathfrak U_c$, $y_*$ is a critical point of
$\Psi_{A,B}$ and its radius is a zero of $g$ in $(r_1,r_2)$.
 
\smallskip\noindent\emph{(b)} Since $g>0$ near $0$ and near $1$, we can choose $r_1<\hat r<r_2$ with
$g(r_1)>0$ and $g(r_2)>0$. We apply (a) on $(r_1,\hat r)$ and on $(\hat r,r_2)$, where the degrees are
\[
  \tfrac12\big(\sign g(\hat r)-\sign g(r_1)\big)=-1,\qquad
  \tfrac12\big(\sign g(r_2)-\sign g(\hat r)\big)=+1.
\]
This gives two solutions $u^\pm_\e$ with radii $r^-_\e\in(r_1,\hat r)$ and $r^+_\e\in(\hat r,r_2)$.
The set of zeros of $g$ in $[r_1,\hat r]$ is compact and contains neither $r_1$ nor $\hat r$, so it lies in
$[r_1,\hat r-\varrho]$ for some $\varrho>0$. Similarly the zeros of $g$ in $[\hat r,r_2]$ lie in $[\hat r+\varrho,r_2]$.
Since $r^\pm_\e$ accumulate only at these zeros, for $\e$ small
\[
  r^-_\e<\hat r-\tfrac\varrho2<\hat r+\tfrac\varrho2<r^+_\e .
\]
The negative parts of $u^-_\e$ and $u^+_\e$ concentrate on rings of radii that stay apart, so the
two solutions are distinct.
 
\smallskip\noindent\emph{(c)} Fix $\eta'\in(0,\frac12)$ with $\mathcal K\subset\interior\mathcal P_{\eta'}$, and argue by
contradiction. Suppose there are $\e_n\to0$ and solutions $u_n=V_{y_n}+\phi_n$ of \eqref{pb} with
$y_n\in \mathcal K$, $\phi_n\in \Zs^\perp$ and $\|\phi_n\|\le C\vartheta(\mu_n)$. Since $u_n$ solves \eqref{pb}, the projected
equation of Proposition~\ref{prop:phi} holds, so $\phi_n=\phi_{y_n}$ by uniqueness. Hence
$\nabla\widetilde J_{\e_n}(y_n)=0$, because $\partial_{y_i}\widetilde J_\e(y)=J_\e'(V_y+\phi_y)[\partial_{y_i}(V_y+\phi_y)]$. By
Proposition~\ref{prop:exp},
\[
  \nabla\Psi_{A,B}(y_n)=-\nabla\varpi_{\e_n}(y_n)\to0 .
\]
Up to a subsequence $y_n\to y_*\in \mathcal K$, and $\nabla\Psi_{A,B}(y_*)=0$. Then the radius of $y_*$ is a zero of $g$,
which contradicts $g>0$.
\end{proof}

\subsubsection{Behaviour of \texorpdfstring{$g$}{g} at the endpoints}\label{ssec:ends}
We now study the behaviour of $g$ at the endpoints of $(0,1)$. Near $r=0$ we use the second expression of $g$ in \eqref{eq:g}, namely
\begin{equation*}
g(r)=1+k\beta(r)\sigma(r)-s(r)^{-2}\sigma(r)\,(\varphi_k(r)\sigma(r)+\beta(r)).
\end{equation*}

\begin{lemma}\label{lem:zero}
For all $N\ge5$ and $k\ge2$, $g(r)\to+\infty$ as $r\to0^+$. More precisely, letting
\[\Lambda_k:=\sum_{j=1}^{k-1}(2\mathfrak s_j)^{-(N-2)}>0 \qquad \text{and} \qquad s_0:=\left(\frac{2}{\Lambda_k}\right)^{\frac2{N-2}},\]
as $r \to 0^{+}$, we have the following expansion
\[g(r)=\frac{2}{\Lambda_k}\big(k+s_0^{-2}\big)\,r^{-(N-2)}+O(1).\]
\end{lemma}
\begin{proof}
As $r \to 0^{+}$ we have $a\to1$ and $b_j\to0$. In \eqref{eq:phiexpl} the only singular terms are
\[b_j^{-(N-2)}=(2\mathfrak s_j)^{-(N-2)}\,r^{-(N-2)};\]
indeed all the other terms are smooth even functions of $r$ near $0$. Hence
\begin{equation}\label{eq:expfi}
\varphi_k=-\Lambda_k\,r^{-(N-2)}+O(1),\qquad \varphi_k'=(N-2)\Lambda_k\,r^{-(N-1)}+O(r).
\end{equation}
By the definition \eqref{eq:s} of $\sigma$ and by \eqref{eq:expfi} we get
\[\sigma=\frac2{\Lambda_k}\big(1+O(r^N)\big),\qquad s=s_0\big(1+O(r^N)\big).\]
Moreover,
\[k\beta(r)\sigma(r) = \frac{2k}{\Lambda_k} \frac{1}{r^{N-2}} + O(1)\quad\text{and}\quad 
\varphi_{k}(r) \sigma(r) = -\frac{2}{r^{N-2}} + O(1);\]
hence
\[\varphi_k(r)\sigma(r)+\beta(r)=-\frac{1}{r^{N-2}}+O(1),\qquad\text{so that}\qquad -s^{-2}\sigma(\varphi_k\sigma+\beta)=\frac{2}{\Lambda_k}\,s_0^{-2}\,\frac1{r^{N-2}}+O(1).\]
Finally, as $r \to 0^{+}$,
\begin{equation*}
g(r) = \frac{2}{\Lambda_{k}} \left(k+\frac{1}{s_0^{2}} \right) \frac{1}{r^{N-2}} + O(1)
\end{equation*}
and the claim follows, since both contributions are positive multiples of $r^{-(N-2)}$.
\end{proof}
Near $r=1$ we use instead the first expression of $g$ in \eqref{eq:g}, namely
\[g(r)=1+k\beta(r) s(r)^{\frac{N-2}{2}}-\varphi_k(r)s(r)^{N-4}-\beta(r) s(r)^{\frac{N-6}{2}}.\]
\begin{lemma}\label{lem:one}
As $r\to1^-$,
\[
g(r)=1-\frac N2\,2^{\mathfrak e_N}(1-r)^{\mathfrak e_N}\big(1+O(1-r)\big)
\]
where
\[\mathfrak e_N:=\frac{N^2-6N+4}{N-2}.\]
In particular, since $\mathfrak e_5=-\frac13<0$ and $\mathfrak e_N>0$ for $N\ge6$, we get, as $r \to 1^{-}$,
\begin{equation*}
g\to-\infty\quad\text{if}\quad N=5\quad\text{and}\quad g\to1\quad\text{if}\quad N\ge6.
\end{equation*}
\end{lemma}

\begin{proof}
As $r\to 1^-$ we have $a=1-r^2 =2(1-r)\big(1+O(1-r)\big)$,
while $b_j\to2\mathfrak s_j>0$. Hence the sums in \eqref{eq:phiexpl} and in \eqref{eq:phiprime} are bounded up to $r=1$, and only the diagonal terms matter:
\[\varphi_k(r)=\frac{1}{2^{N-2}(1-r)^{N-2}}\big(1+O(1-r)\big),\qquad \varphi'_{k}(r)=\frac{N-2}{2^{N-2}(1-r)^{N-1}}\big(1+O(1-r)\big),\]
and
\[\beta(r)=(N-2)(1-r)+O\big((1-r)^2\big),\qquad |\beta'(r)|=(N-2)+O(1-r).\]
Consequently
\[
\sigma=\frac{2|\beta'|}{\varphi_k'}=2^{N-1}(1-r)^{N-1}\left(1+O\left((1-r)\right)\right).
\]
Since $s=\sigma^{\frac2{N-2}}$, a direct computation gives
\[k\beta s^{\frac{N-2}2}=k\beta\sigma=O\big((1-r)^N\big),\qquad \varphi_ks^{N-4}=2^{\mathfrak e_N}(1-r)^{\mathfrak e_N}\big(1+O(1-r)\big)\]
and
\[\beta s^{\frac{N-6}2}=\frac{N-2}2\, 2^{\mathfrak e_N}(1-r)^{\mathfrak e_N}\big(1+O(1-r)\big),\]
because $\frac{2(N-1)(N-4)}{N-2}-(N-2)=1+\frac{(N-1)(N-6)}{N-2}=\mathfrak e_N$. Inserting these expansions in the expression of $g$, the claim follows.
\end{proof}

\begin{remark}
By Lemmas~\ref{lem:zero} and~\ref{lem:one}, in dimension $N=5$ the function $g$ changes sign on $(0,1)$ for every $k$; this is the source of Theorem~\ref{thm:main1}. For $N\ge6$, instead, $g$ is positive near both endpoints, so a sign change can only come from the interior of $(0,1)$.\end{remark}
\section{Proof of Theorem~\ref{thm:main1}\,(i)}\label{sec:N5}

Let $N=5$. By Lemmas~\ref{lem:zero} and~\ref{lem:one},
\[g(r)\to +\infty\quad\text{as}\quad r\to 0^+,\qquad g(r)\to -\infty\quad\text{as}\quad r\to 1^-.\]
Choose $0<r_1<r_2<1$ with $g>0$ on $(0,r_1]$ and $g<0$ on $[r_2,1)$.
By Proposition~\ref{prop:degree},
\[
\deg(\nabla\Psi,\mathfrak U,0)=\tfrac12\big(\sign g(r_2)-\sign g(r_1)\big)=\tfrac12(-1-1)=-1\neq0,
\]
and Corollary~\ref{cor:exist}\,(a) gives the solution. Since $\mu=\e$ when $N=5$, the estimate $\|\phi_\e\|_{H^1_0}\le C\vartheta(\mu)=C\e^{5/2}$ follows from Proposition~\ref{prop:phi}. Part (ii) of Theorem~\ref{thm:main1} will be proved in Section~\ref{sec:lower} (Proposition~\ref{prop:N5ring}).

\begin{remark}\label{rem:morse}
Let $r_*$ be a zero of $g$ with $g'(r_*)<0$ (by Lemmas~\ref{lem:zero} and~\ref{lem:one}, if all the zeros of $g$ are simple, at least one of them has this property). Then the corresponding critical point $y_*$ of $\Psi$ is a non-degenerate saddle point with Morse index one. Indeed, $\det D^2\Psi(y_*)<0$ by Lemma~\ref{lem:crit}\,(ii), so $D^2\Psi(y_*)$ has one or three negative eigenvalues. It cannot have three, because $\Psi$ is convex along the dilation $t\mapsto(r,td_0,td_1)$ at $t=1$: writing $\Psi(r,td_0,td_1)=t^{N-2}E_1-t^2E_2$ with $E_2=\frac12(d_0^2+kd_1^2)>0$, criticality gives $(N-2)E_1=2E_2$, and then
\[
\frac{d^2}{dt^2}\Big|_{t=1}\Psi(r,td_0,td_1)=(N-2)(N-3)E_1-2E_2=2(N-4)\,E_2>0 .
\]
\end{remark}

\begin{remark}\label{rem:N5num}
Numerically, for $N=5$ the function $g$ has exactly one zero $r_k$ for every $k\ge2$, and $g'(r_k)<0$; hence, by Remark~\ref{rem:morse}, the solution of Theorem~\ref{thm:main1} comes from a non-degenerate critical point of Morse index one, and $r_\e\to r_k$. 
\end{remark}

\section{Proof of Theorem~\ref{thm:main2}\,(i)}\label{sec:k6}

Let $N\ge6$. By Lemmas~\ref{lem:zero} and~\ref{lem:one},
\[g(r)\to +\infty\quad\text{as}\quad r\to 0^+,\qquad g(r)\to 1\quad\text{as}\quad r\to 1^-.\]
Hence, by Corollary~\ref{cor:exist}\,(b), there are two solutions as soon as $g(\hat r)<0$ at a single point $\hat r$. We choose
\[
\hat r=2^{-1/2},\qquad\text{i.e.}\qquad a(\hat r)=1-\hat r^2=\tfrac12=\hat r^2
\]
since at this radius the interaction with the centre (measured by $\beta$) and the interaction with the boundary (measured by $a^{-(N-2)}$) are of the same size, and all the quantities appearing in $g$ become explicit powers of $2$. Theorem~\ref{thm:main2}\,(i) then follows from Proposition~\ref{prop:k6} below.

\subsection{An explicit formula for $g(2^{-1/2})$}

Throughout the rest of the paper we use the notation
\begin{equation}\label{eq:mw}
m:=\frac{N-2}{2},\qquad \ell:=\frac{4}{N-2}=\frac2m ,
\end{equation}
so that $m\ge2$ and $\ell\in(0,1]$ when $N\ge6$. In this section we also set
\[
\chi:=2^{-m},\qquad u_j:=8\mathfrak s_j^2,\qquad j=1,\dots,k-1,
\]
and
\begin{equation}\label{eq:niD}
\mathfrak n_{k,m}:=\sum_{j=1}^{k-1}\Big[u_j^{-m}-(1+u_j)^{-m}\Big],\qquad
D:=1+\sum_{j=1}^{k-1}\Big[\frac12\,u_j^{-m}-\Big(\frac{u_j}{2}-1\Big)(1+u_j)^{-(m+1)}\Big].
\end{equation}
Clearly $\mathfrak n_{k,m}\ge 0$. Moreover $D\ge 1$, since (see the proof below) the $j$-th term of
$D-1$ is $\omega_j$ divided by the diagonal term of $\varphi_k'$ at $r=2^{-1/2}$, and $\omega_j>0$ by
Lemma~\ref{lem:phiprime}.

\begin{lemma}\label{lem:formula}
We have
\begin{equation}\label{eq:gval}
g(2^{-1/2})=1-4(1-\chi)\,D^{\ell-1}+\frac1D\Big[k(1-\chi)-4(1-\mathfrak n_{k,m})\,D^{\ell-1}\Big].
\end{equation}
\end{lemma}

\begin{proof}
At $\hat r=2^{-1/2}$ we have $a=\hat r^{2}=\frac12$, $b_j^2=4\hat r^2\mathfrak s_j^2=\frac{u_j}{4}$
and $a^2+b_j^2=\frac{1+u_j}{4}$. Hence, by \eqref{eq:phiexpl}, \eqref{eq:signs} and
Lemma~\ref{lem:phiprime}, all the powers of $4$ factor out and
\[
\varphi_k=4^m(1-\mathfrak n_{k,m}),\qquad \beta=2^m-1=\frac{1-\chi}{\chi},\qquad
|\beta'|=2m\sqrt2\,2^{m},\qquad \varphi_k'=4m\sqrt2\,4^{m}\,D .
\]
For $\varphi_k'$: the diagonal term is $2(N-2)\hat r a^{-(N-1)}=4m\sqrt2\,4^m$, and a direct
computation gives
\[
\frac{\omega_j}{4m\sqrt2\,4^m}=\frac12\,u_j^{-m}-\Big(\frac{u_j}{2}-1\Big)(1+u_j)^{-(m+1)} .
\]
Therefore
\[
\sigma=\frac{2|\beta'|}{\varphi_k'}=\frac{\chi}{D},\qquad s^{-2}=\sigma^{-\ell}=4D^{\ell}.
\]
Since $4^m \chi^2=1$, we get $\varphi_k\sigma^2=\frac{1-\mathfrak n_{k,m}}{D^2}$ and $\beta\sigma=\frac{1-\chi}{D}$.
Inserting these identities in $g=1+k\beta\sigma-s^{-2}\sigma(\varphi_k\sigma+\beta)$, see \eqref{eq:g}, the claim follows.
\end{proof}

\begin{proposition}\label{prop:k6}
$g(2^{-1/2})<0$ for every $N\ge 6$ and $2\le k\le 6$.
\end{proposition}

\begin{proof}
If $2\le k\le 6$, then for $1\le j\le k-1$ we have
$\mathfrak s_j\ge\sin\frac{\pi}{k}\ge\frac12$, i.e.\ $u_j\ge 2$, so that
$\frac{u_j}{2}-1\ge 0$. Letting $\Sigma_{k,m}:=\sum_{j=1}^{k-1}u_j^{-m}$, by \eqref{eq:niD} we obtain
\begin{equation}\label{eq:DnuS}
1\le D\le 1+\frac{\Sigma_{k,m}}{2},\qquad 0\le \mathfrak n_{k,m}\le \Sigma_{k,m}.
\end{equation}
Since $u_j\ge2$ and $u\mapsto u^{-m}$, $u\mapsto u^{-m}-(1+u)^{-m}$ are decreasing on $[2,\infty)$, one checks that, for $2\le k\le6$, both $\Sigma_{k,m}$ and $\mathfrak n_{k,m}$ are maximal for $k=6$; moreover they are decreasing in $m$. Hence
\[
\begin{aligned}
\Sigma_{k,m}&\le \Sigma_{6,m}=2\cdot 2^{-m}+2\cdot 6^{-m}+8^{-m},\\
\mathfrak n_{k,m}&\le \mathfrak n_{6,m}=2\big(2^{-m}-3^{-m}\big)+2\big(6^{-m}-7^{-m}\big)+\big(8^{-m}-9^{-m}\big).
\end{aligned}
\]
Moreover, by \eqref{eq:DnuS},
\begin{equation}\label{eq:stimaD}
D^{\ell-1}\ge \Big(1+\frac {\Sigma_{k,m}}{2}\Big)^{\ell-1}\ge 1-\frac{1-\ell}2\, \Sigma_{k,m} .
\end{equation}

\smallskip\noindent\emph{The case $N=6$.} Here $m=2$, $\ell=1$ and $\chi=\frac14$, so that \eqref{eq:gval} gives
\[
g(2^{-1/2})=-2+\frac{\frac 34 k-4+4\mathfrak n_{k,2}}{D}<-2+\frac 2 D\le 0,
\]
because $k\le 6$, $\mathfrak n_{k,2}\le\mathfrak n_{6,2}=\frac{75137}{254016}<\frac38$ and $D\ge 1$.

\smallskip\noindent\emph{The cases $N=7,8,9$.} The quantities $m$, $\ell$ and $\chi$ are explicit, and a computer-assisted evaluation of \eqref{eq:gval} gives $g(2^{-1/2})<0$.

\smallskip\noindent\emph{The case $N\ge 10$.} Here $m\ge4$. By \eqref{eq:gval}, $D\ge1$ and \eqref{eq:stimaD},
\[
g(2^{-1/2})\le 1-4(1-\chi)\,D^{\ell-1}+\Big[k(1-\chi)-4(1-\mathfrak n_{k,m})\,D^{\ell-1}\Big]^+ .
\]
If the bracket is nonpositive, then $g(2^{-1/2})\le1-4(1-\chi)\big(1-\frac12\Sigma_{6,4}\big)\le1-\frac{15}4\cdot\frac34<0$, since $\Sigma_{6,4}<\frac12$. Otherwise, using $k\le6$ and \eqref{eq:stimaD},
\[
\begin{aligned}
g(2^{-1/2})&\le 1-4(1-\chi)\Big(1-\frac{1-\ell}2\, \Sigma_{6,m}\Big)+6(1-\chi)-4(1-\mathfrak n_{6,m})\Big(1-\frac{1-\ell}2\, \Sigma_{6,m}\Big)\\
&\le -1-2\chi+4\mathfrak n_{6,m}+4(1-\ell)\,\Sigma_{6,m}
\ <\ -1+4\mathfrak n_{6,4}+4\Sigma_{6,4}<0,
\end{aligned}
\]
because $\mathfrak n_{6,m}$ and $\Sigma_{6,m}$ are decreasing in $m$, and $4\mathfrak n_{6,4}+4\Sigma_{6,4}<0.92$.
\end{proof}

\section{Lower bounds for \texorpdfstring{$g$}{g}: non-existence for large \texorpdfstring{$k$}{k} and localisation for \texorpdfstring{$N=5$}{N=5}}\label{sec:lower}
We first rewrite $g$ in a form which is more suitable to obtain lower bounds. Let
\[\mathcal Q(r) := 2 |\beta'(r)| \varphi_{k}(r) + \beta(r) \varphi'_{k}(r) \qquad \text{and} \qquad \gamma(r):=\frac{\beta(r)}{|\beta'(r)|} = \frac{r(1-r^{N-2})}{N-2}.\]
Since $1-r^{N-2} \le (N-2)(1-r)$, we have
\[\gamma(r) \le \frac{r}{N-2} \qquad \text{and} \qquad \gamma(r) \le r(1-r) \le a.\]
Furthermore, by \eqref{eq:s},
\[\varphi_{k}(r)\sigma(r)+\beta(r) = \frac{\mathcal Q(r)}{\varphi'_{k}(r)} \qquad \text{and} \qquad s(r)^{-2} \sigma(r) = \sigma(r)^{1-\ell},\]
with $\ell$ as in \eqref{eq:mw}, so that \eqref{eq:g} becomes
\begin{equation}\label{eq:gN}
g(r) = 1+k\beta(r) \sigma(r) - \sigma(r)^{1-\ell} \frac{\mathcal Q(r)}{\varphi'_{k}(r)}.
\end{equation}
Since $k\beta\sigma>0$ and $\sigma, \varphi'_{k}>0$, we get the following criterion, valid for every $N\ge5$:
\begin{equation}\label{eq:Ncrit}
\mathcal Q(r) \le 0 \quad\Longrightarrow\quad g(r) \ge 1.
\end{equation}
Recalling the definitions of $\tau_{j}$ and $\omega_{j}$ in \eqref{eq:phiexpl} and \eqref{eq:omegaj}, we get
\begin{equation}\label{eq:Nsplit}
\frac{\mathcal Q(r)}{|\beta'(r)|}=2\varphi_k(r)+\gamma\varphi_k'(r)=\Big[2a^{-(N-2)}+2(N-2)\,\gamma\,r\,a^{-(N-1)}\Big]-\sum_{j=1}^{k-1}\big(2\tau_j-\gamma\,\omega_j\big).
\end{equation}

\subsection{The non-existence theorem}

\begin{lemma}\label{lem:bulk}
Let $N\ge5$, $k\ge2$,
\[c_N:=\int_0^\infty(1+z^2)^{-\frac{N-2}2}\dd z = \frac{\sqrt{\pi}}{2}\, \frac{\Gamma\left(\frac{N-3}{2}\right)}{\Gamma \left(\frac{N-2}{2}\right)}\]
and
\[F_N(M):=2(N-1)\left(\frac{2\pi}M\right)^{N-2}+(N-1)\,c_N\,(2\pi)^{N-2}M^{-(N-3)}, \qquad M>0.\]
If $F_N(M)<1$, then $\mathcal Q(r)<0$, and hence $g(r)\ge1$, whenever $k(1-r^2)\ge M$.
\end{lemma}

\begin{proof}
We estimate the two parts of \eqref{eq:Nsplit}. First, since $\gamma\le a$ and $r\le1$,
\begin{equation}\label{eq:aa}
2a^{-(N-2)}+2(N-2)\,\gamma\,r\,a^{-(N-1)} \le 2a^{-(N-2)}+2(N-2)a^{-(N-2)}=2(N-1)a^{-(N-2)}.
\end{equation}
Next, by the definition \eqref{eq:omegaj} of $\omega_{j}$ and since $\gamma\le\frac r{N-2}$, $\gamma\le a$, $r\le1$ and $a^2\le a^2+b_j^2$, we get
\[\gamma\,\omega_j\le b_j^{-(N-2)}-b_j^2(a^2+b_j^2)^{-\frac N2}+2(N-2)(a^2+b_j^2)^{-\frac{N-2}2}.\]
Since $2\tau_j=2b_j^{-(N-2)}-2(a^2+b_j^2)^{-\frac{N-2}2}$, it follows that
\begin{equation}\label{eq:tauomega}
2\tau_j-\gamma\omega_j\ge b_j^{-(N-2)}-2(N-1)(a^2+b_j^2)^{-\frac{N-2}2}.
\end{equation}
Therefore, by \eqref{eq:aa} and \eqref{eq:tauomega},
\begin{equation}\label{eq:first}
\frac{\mathcal Q(r)}{|\beta'(r)|}\le2(N-1)\,a^{-(N-2)}\;-\;\sum_{j=1}^{k-1}b_j^{-(N-2)}\;+\;2(N-1)\sum_{j=1}^{k-1}(a^2+b_j^2)^{-\frac{N-2}2}.
\end{equation}
We estimate the two sums in \eqref{eq:first}. Since $\sin x\le x$,
\[ \sum_{j=1}^{k-1} b_{j}^{-(N-2)} \ge b_{1}^{-(N-2)} \ge \left( \frac{k}{2\pi r} \right)^{N-2}.\]
Moreover, $b_{k-j}=b_j$ and hence
\[\sum_{j=1}^{k-1}(a^2+b_j^2)^{-\frac{N-2}2}\le2\sum_{j=1}^{\lfloor k/2\rfloor} (a^2+b_j^2)^{-\frac{N-2}2}.\]
For $1\le j\le\lfloor\frac k2\rfloor$ we have $\frac{\pi j}k\in(0,\frac\pi2]$, so that, by Jordan's inequality $\sin x \ge \frac{2}{\pi}x$ on $[0,\frac\pi2]$, $\mathfrak s_j\ge\frac{2j}k$ and $b_j\ge\frac{4rj}k$. Hence, letting
\[f(\zeta):=\Big(a^2+\frac{16r^2\zeta^2}{k^2}\Big)^{-\frac{N-2}2},\qquad \zeta\ge0,\]
we have $(a^2+b_j^2)^{-\frac{N-2}2}\le f(j)$. Since $f$ is decreasing on $[0,\infty)$, $f(j)\le\int_{j-1}^jf(\zeta)\dd \zeta$ for $j\ge1$, and
\[\sum_{j=1}^{\lfloor k/2\rfloor} f(j)\le\int_0^\infty f(\zeta) \dd \zeta .\]
With the change of variables $\zeta=\frac{ka}{4r}u$ we have $a^2+\frac{16r^2\zeta^2}{k^2}=a^2(1+u^2)$, so that
\[\int_0^\infty f(\zeta)\dd \zeta=a^{-(N-2)}\,\frac{ka}{4r}\int_0^\infty(1+u^2)^{-\frac{N-2}2}\dd u=\frac{k\,c_N}{4r}\,a^{-(N-3)}.\]
Putting everything together,
\[
\sum_{j=1}^{k-1}(a^2+b_j^2)^{-\frac{N-2}2}\le\frac{k\,c_N}{2r}\,a^{-(N-3)} .
\]
Finally, multiplying \eqref{eq:first} by the positive number $\big(\frac{2\pi r}k\big)^{N-2}$ and using $r\le1$, we obtain
\begin{equation*}
\left(\frac{2\pi r}{k} \right)^{N-2} \frac{\mathcal Q(r)}{|\beta'(r)|} \le \frac{2(N-1)(2\pi)^{N-2}}{(ka)^{N-2}} -1 + (N-1)c_{N} (2\pi)^{N-2} (ka)^{-(N-3)} = F_{N}(ka)-1.
\end{equation*}
The function $F_N$ is strictly decreasing on $(0,\infty)$, because both its terms are negative powers of $M$. Hence, if $ka\ge M$, then $F_N(ka)\le F_N(M)<1$, the right-hand side above is negative and $\mathcal Q(r)<0$. The last claim follows from \eqref{eq:Ncrit}.
\end{proof}

\begin{lemma}\label{lem:layer}
Let $N\ge6$. For all $r\in(0,1)$ and $k\ge2$,
\[
g(r)\ge1-(N-1)\,r^{-\frac{2N(N-4)}{N-2}}\,a^{\mathfrak e_N},\qquad \mathfrak e_N=\frac{N^2-6N+4}{N-2}.
\]
\end{lemma}

\begin{proof}
By \eqref{eq:gN}, the first expression in \eqref{eq:g} and $k\beta\sigma\ge0$, we get
\begin{equation}\label{eq:glow}
g\ge1-\sigma^{1-\ell}\big(\varphi_k\sigma+\beta\big).
\end{equation}
If $\varphi_k(r)\sigma(r)+\beta(r)\le0$ the claim is obvious, so we assume $\varphi_k(r)\sigma(r)+\beta(r)>0$.
Let $\sigma_0(r):=r^{-N}a^{N-1}$. By \eqref{eq:phiprime}, $\varphi'_k(r) \ge 2(N-2)r\,a^{-(N-1)}$, hence $\sigma(r)\le\sigma_0(r)$.
Moreover $\varphi_k(r)\le a^{-(N-2)}$ by \eqref{eq:phiexpl}, and
\begin{equation}\label{eq:betaa}
\beta=r^{-(N-2)}(1-r^{N-2})\le\frac{N-2}2\,a\,r^{-(N-2)},
\end{equation}
because $1-y^{q}\le q(1-y)$ for all $y \in[0,1]$ and $q \ge 1$; it suffices to take $y=r^2$ and $q=\frac{N-2}2$. Consequently
\begin{equation}\label{eq:phikrb}
\varphi_k\sigma+\beta\le a^{-(N-2)}\sigma_0+\tfrac{N-2}2\,a\,r^{-(N-2)}= r^{-N}a+\tfrac{N-2}2\,a\,r^{-(N-2)}\le(N-1)\,r^{-N}a
\end{equation}
and, since $1-\ell\ge0$,
\begin{equation}\label{eq:sigma1r}
\sigma^{1-\ell}\le\sigma_0^{1-\ell}=r^{-\frac{N(N-6)}{N-2}}a^{\frac{(N-1)(N-6)}{N-2}}.
\end{equation}
Combining \eqref{eq:glow}, \eqref{eq:phikrb} and \eqref{eq:sigma1r}, the claim follows.
\end{proof}

\begin{theorem}\label{thm:KN}
Let $N\ge6$. Let $M>0$ satisfy $F_N(M)<1$, and let $K>M$ satisfy
\[
L_N(M,K):=(N-1)\Big(1-\frac MK\Big)^{-\frac{N(N-4)}{N-2}}\Big(\frac MK\Big)^{\mathfrak e_N}<1 .
\]
Then $g>0$ on $(0,1)$ for every $k\ge K$. Such $M$ and $K$ exist for every $N\ge6$; we denote by $K_N$ the smallest integer $K$ obtained in this way.
\end{theorem}

\begin{proof}
Fix $r\in(0,1)$ and $k\ge K$. If $ka\ge M$, then $g(r)\ge1$ by Lemma~\ref{lem:bulk}. Otherwise $a<M/k\le M/K$; then $r^2=1-a>1-\frac MK$ and, since $\mathfrak e_N>0$, Lemma~\ref{lem:layer} gives $g(r)\ge1-L_N(M,K)>0$.
As for the existence of $M$ and $K$: $F_N(M)\to0$ as $M\to\infty$, so there exists $\bar{M}$ with $F_N(\bar M)<1$, and $L_N(\bar{M},K)\to0$ as $K\to\infty$, because $\mathfrak e_N>0$.
\end{proof}

\begin{proof}[Proof of Theorem~\ref{thm:main2}\,\textup{(ii)}]
By Theorem~\ref{thm:KN}, $g>0$ on $(0,1)$ for every $k\ge K_N$, and Corollary~\ref{cor:exist}\,(c) gives the claim. Since by Proposition~\ref{prop:k6} $g$ changes sign for $2\le k\le 6$, necessarily $K_N\ge7$.
\end{proof}

\subsection{Dimension five: the ring approaches the boundary}
In dimension $N=5$ Lemma~\ref{lem:bulk} still holds and gives an explicit localisation of the zeros of $g$.

\begin{proposition}\label{prop:N5ring}
Let $N=5$ and $k\ge2$. Then $g(r)\ge1$ whenever $k(1-r^2)\ge 33$. In particular every zero $r_k$ of $g$ satisfies
\[
1-r_k^2<\frac{33}{k},
\]
and Theorem~\ref{thm:main1}\,\textup{(ii)} holds.
\end{proposition}

\begin{proof}
For $N=5$ we have $c_5=\int_0^\infty(1+z^2)^{-3/2}\dd z=1$, so that
\[
F_5(M)=\frac{8(2\pi)^3}{M^3}+\frac{4(2\pi)^3}{M^2}.
\]
Since $(2\pi)^3<249$, we get $F_5(33)<1$, and Lemma~\ref{lem:bulk} with $M=33$ gives the first claim. Hence every zero of $g$ satisfies $k(1-r^2)<33$.

Now let $u_\e$ be the solutions of Theorem~\ref{thm:main1}\,(i), with radii $r_\e$. By Corollary~\ref{cor:exist}\,(a), every limit point $r_*$ of $r_\e$ as $\e\to0$ is a zero of $g$, so that $k(1-r_*^2)<33$. If $k(1-r_{\e_n}^2)\ge33$ along a sequence $\e_n\to0$, then, up to a subsequence, $r_{\e_n}\to r_*$ with $k(1-r_*^2)\ge33$, a contradiction. Thus $k(1-r_\e^2)<33$ for $\e$ small, and every limit point $r_k$ of $r_\e$ satisfies $1-r_k^2<33/k$, which tends to $0$ as $k\to+\infty$.
\end{proof}

\begin{remark}
The constant $33$ is far from optimal: by Remark~\ref{rem:N5num}, numerically $k(1-r_k^2)\to 8$ approximately. The proof shows the mechanism behind Theorem~\ref{thm:main1}\,(ii): when the distance $a\approx2(1-r)$ of the ring from the boundary is large compared with the distance $\approx 2\pi r/k$ between consecutive negative bubbles, the mutual interaction of the negative bubbles (the ring part of $\varphi_k$ in \eqref{eq:Nsplit}) dominates their interaction with the boundary, and the reduced function has no critical points. The ring can only find an equilibrium where these two interactions are comparable, that is, at distance $O(1/k)$ from $\partial\mathcal B_1$.
\end{remark}

\begin{figure}[htbp]
    \centering
    \includegraphics[width=0.5\textwidth]{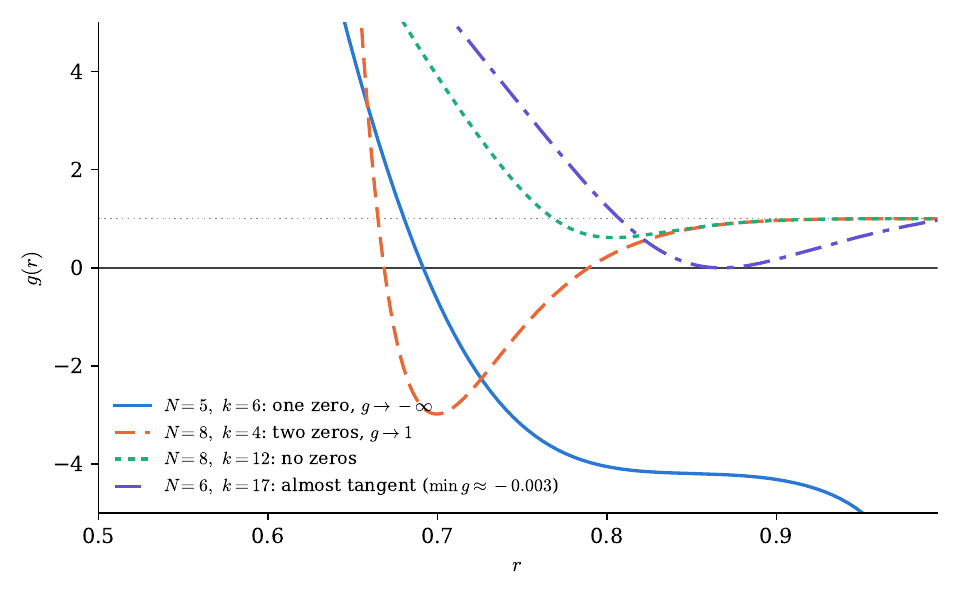}
    \caption{\tiny{The possible behaviours of the function $g$.}}
    \label{fig:g}
\end{figure}

The possible behaviours of $g$ described in Sections~\ref{sec:Psi}--\ref{sec:lower} are summarised in Figure~\ref{fig:g}.

\section{The nodal set near the boundary}\label{sec:nodal}

In this section we prove Theorem~\ref{thm:nodal}. We use the notation of the previous sections and, in addition, we set
\[
t:=r^2,\qquad \theta_j:=\frac{2\pi j}{k},\qquad e_j:=(\cos\theta_j,\sin\theta_j,0,\dots,0),\qquad q_j:=|1-t\,e^{i\theta_j}|^2,\qquad 0\le j\le k-1,
\]
so that $\xi_{j+1}(r)=r e_j$, $q_0=a^2$ and $q_j=1-2t\cos\theta_j+t^2=a^2+b_j^2$ for $1\le j\le k-1$. Note that $\{e_0,\dots,e_{k-1}\}=\{\xi_j/|\xi_j|:\ j=1,\dots,k\}$. We denote by
$\nu$ the exterior unit normal to $\partial\mathcal B_1$, and we write the points of $\partial\mathcal B_1$ as $x=(z,x'')$, with $z\in\R^2$, $|z|\le1$.

\subsection{The normal derivative on \texorpdfstring{$\partial \mathcal B_1$}{the boundary}}
For $r\in(0,1)$ and $x\in\partial\mathcal B_1$ we define the \emph{boundary density}
\begin{equation}\label{eq:P}
P_r(x):=(1-r^2)\sum_{j=1}^k|x-\xi_j(r)|^{-N}.
\end{equation}
Recall that for $|x|=1$ and $|\xi|<1$ one has
$\partial_\nu G(x,\xi)=-(N-2)(1-|\xi|^2)|x-\xi|^{-N}$ (Poisson kernel).

The following result is standard: away from the concentration points, $\mu^{-\frac{N-2}2}u_\e$ converges in $C^1$ up to the boundary to the corresponding combination of Green functions. It follows from the expansion $PU_{\delta,\xi}=\alpha_N\delta^{\frac{N-2}2}G(\cdot,\xi)+O(\delta^{\frac{N+2}2})$ in $C^1$ away from $\xi$, from the estimate $\|\phi_\e\|_{H^1_0}=o(\mu^{\frac{N-2}2})$ of Proposition~\ref{prop:phi}, and from elliptic regularity near $\partial\mathcal B_1$ (see e.g.\ \cite{BDP, R, MP}).

\begin{proposition}\label{prop:normal}
Let $\e_n\to0$ and let $y_n=(r_n,d_{0,n},d_{1,n})\in\mathcal P_\eta$ be such that $u_n:=V_{y_n}+\phi_{y_n}$ solves \eqref{pb} with $\e=\e_n$. Assume $y_n\to y_*=(r_*,d_{0*},d_{1*})$. Then, uniformly on $\partial\mathcal B_1$,
\[
\mu_n^{-\frac{N-2}{2}}\,\partial_\nu u_n\ \longrightarrow\ (N-2)\,\alpha_N\,d_{0*}^{\frac{N-2}{2}}
\Big[\Big(\frac{d_{1*}}{d_{0*}}\Big)^{\frac{N-2}{2}}P_{r_*}(x)-1\Big].
\]
If $y_*$ is a critical point of $\Psi_{A,B}$, then $(d_{1*}/d_{0*})^{\frac{N-2}{2}}=\sigma(r_*)$ and
$g(r_*)=0$.
\end{proposition}

\begin{proof}
As recalled above, $\mu_n^{-\frac{N-2}{2}}u_n\to\alpha_N\big[d_{0*}^{\frac{N-2}{2}}G(\cdot,0)-d_{1*}^{\frac{N-2}{2}}
\sum_jG(\cdot,\xi_j(r_*))\big]$ in $C^1$ of a neighbourhood of $\partial\mathcal B_1$, and the formula for $\partial_\nu G$ gives the limit. If $y_*$ is critical for $\Psi_{A,B}$, then
$(r_*,d_{0*}/c,d_{1*}/c)$ is critical for $\Psi$ (see Section~\ref{sec:prelim}), so by Lemma~\ref{lem:crit} $g(r_*)=0$ and
$d_{1*}/d_{0*}=s(r_*)$, i.e.\ $(d_{1*}/d_{0*})^{\frac{N-2}{2}}=s(r_*)^{\frac{N-2}{2}}=\sigma(r_*)$ by \eqref{eq:s}.
\end{proof}

Thus, near $\partial\mathcal B_1$, the sign of $\partial_\nu u_\e$ is governed by the sign of
$\sigma(r_*)P_{r_*}-1$: where $\sigma P<1$ we have $\partial_\nu u_\e<0$, i.e.\ $u_\e>0$ near that
portion of $\partial\mathcal B_1$, and where $\sigma P>1$ we have $u_\e<0$ near $\partial\mathcal B_1$.

\subsection{Minimum and maximum of \texorpdfstring{$P_r$}{P}}

Note that $P_r(x)=-\frac1{N-2}\sum_{j}\partial_\nu G(x,\xi_j(r))$ is, up to the factor $|\partial\mathcal B_1|$, a sum of $k$ Poisson kernels of the ball. For $\rho\in(0,1)$ we set
\begin{equation}\label{eq:S}
\mathcal S(\rho):=\sum_{j=0}^{k-1}|1-\rho\,e^{i\theta_j}|^{-N},\qquad\text{so that}\qquad \sum_{j=0}^{k-1}q_j^{-\frac N2}=\mathcal S(t),
\end{equation}
and
\[
c_n:=\frac{(N/2)_n}{n!}>0,\qquad\text{so that}\qquad (1-\zeta)^{-\frac N2}=\sum_{n\ge0}c_n\zeta^n\quad\text{for }|\zeta|<1 .
\]
The extremal values of $P_r$ are described by the following lemma.

\begin{lemma}\label{lem:extrema}
Let $k\ge2$ and $r\in(0,1)$. Then
\[
\min_{\partial\mathcal B_1}P_r=k\,a\,(1+r^2)^{-\frac N2},\qquad
\max_{\partial\mathcal B_1}P_r=a\,\mathcal S(r),
\]
the minimum being attained on $\partial\mathcal B_1\cap\{x_1=x_2=0\}$ and the maximum exactly at the points $e_0,\dots,e_{k-1}$. Moreover, $\frac1{|\partial\mathcal B_1|}\int_{\partial\mathcal B_1}P_r\,dS=k$.
\end{lemma}

\begin{proof}
\emph{Reduction to the unit disc.} Write $x=(z,x'')\in\partial\mathcal B_1$ with $z\in\overline{\mathbb D}$, where $\mathbb D$ is the unit disc of $\R^2$. Since $|x|=1$ and $e_j$ lies in the $(x_1,x_2)$-plane, $|x-re_j|^2=1+r^2-2r\,z\cdot e_j$, so that $P_r(x)=a\,F(z)$, where
\[
F(z):=\sum_{j=0}^{k-1}f(z\cdot e_j),\qquad f(\lambda):=(1+r^2-2r\lambda)^{-\frac N2},\quad \lambda\in[-1,1].
\]
The function $f$ is positive, increasing and convex, so $F$ is convex on $\overline{\mathbb D}$. Moreover, if $R$ denotes the rotation of angle $\frac{2\pi}k$ in $\R^2$, then $R$ permutes the vectors $e_j$, hence $F\circ R=F$.

\emph{Minimum.} For every $z\in\overline{\mathbb D}$ the points $z,Rz,\dots,R^{k-1}z$ are the vertices of a regular polygon centred at $0$ (or all coincide with $0$), so $\frac1k\sum_{l=0}^{k-1}R^lz=0$. By the invariance and the convexity of $F$,
\[
F(z)=\frac1k\sum_{l=0}^{k-1}F(R^lz)\ \ge\ F\Big(\frac1k\sum_{l=0}^{k-1}R^lz\Big)=F(0)=k\,(1+r^2)^{-\frac N2}.
\]
The value $F(0)$ is attained at $z=0$, that is, on $\partial\mathcal B_1\cap\{x_1=x_2=0\}$.

\emph{Maximum on $\partial\mathbb D$.} For $z=e^{i\theta}$ we have $z\cdot e_j=\cos(\theta-\theta_j)$, hence
\[
F(e^{i\theta})=\sum_{j=0}^{k-1}|1-re^{i(\theta-\theta_j)}|^{-N}.
\]
Since $|1-\zeta|^{-N}=\big|(1-\zeta)^{-\frac N2}\big|^2$, multiplying the series of $(1-\zeta)^{-\frac N2}$ by its complex conjugate we get, for $\alpha\in\R$,
\[
|1-re^{i\alpha}|^{-N}=\sum_{n\ge0}\sum_{p\in\mathbb Z}c_nc_{n+|p|}\,r^{2n+|p|}\,e^{ip\alpha}.
\]
As $\sum_{j=0}^{k-1}e^{-ip\theta_j}$ equals $k$ if $k$ divides $p$ and $0$ otherwise, only the frequencies $p=lk$, $l\in\mathbb Z$, survive:
\begin{equation}\label{eq:Fcos}
F(e^{i\theta})=k\sum_{n\ge0}\sum_{l\in\mathbb Z}c_nc_{n+|l|k}\,r^{2n+|l|k}e^{ilk\theta}
=k\sum_{n\ge0}c_n^2r^{2n}+2k\sum_{l\ge1}\sum_{n\ge0}c_nc_{n+lk}\,r^{2n+lk}\cos(lk\theta).
\end{equation}
All the coefficients are positive, so $F(e^{i\theta})\le F(1)$, with equality if and only if $\cos(k\theta)=1$, i.e.\ $\theta\in\{\theta_0,\dots,\theta_{k-1}\}$. In particular $F$ is not constant on $\partial\mathbb D$, so that $F(1)>F(0)=\min_{\overline{\mathbb D}}F$.

\emph{Maximum on $\overline{\mathbb D}$.} Let $z\in\mathbb D$ and write $z=\rho_0e^{i\theta}$ with $|\rho_0|<1$. The function $\rho\mapsto F(\rho e^{i\theta})$ is convex on $[-1,1]$; if it attained its maximum over $\overline{\mathbb D}$ at the interior point $\rho_0$, it would be constant on $[-1,1]$, and then $F(0)=F(z)=F(1)$, which is impossible. Hence $\max_{\overline{\mathbb D}}F=F(1)$ is attained exactly at the points $z=e^{i\theta_j}$, which correspond to $x=e_j$ (since $|z|=1$ forces $x''=0$). Finally $F(1)=\mathcal S(r)$.

\emph{Mean value.} For every $\xi\in\mathcal B_1$, $\frac{1-|\xi|^2}{|\partial\mathcal B_1|\,|x-\xi|^N}$ is the Poisson kernel of $\mathcal B_1$, whose integral over $\partial\mathcal B_1$ is $1$ (it represents the harmonic function $\equiv1$). Applying this to each $\xi=re_j$ gives the claim.
\end{proof}

\begin{remark}
For $k\ge3$ the vectors $e_j$ span $\R^2$ and $F$ is strictly convex, so the minimum is attained only on $\partial\mathcal B_1\cap\{x_1=x_2=0\}$. For $k=2$, instead, $F(z)$ depends only on $z\cdot e_0=x_1$ and the minimum is attained on the larger sphere $\partial\mathcal B_1\cap\{x_1=0\}$.
\end{remark}

We set
\begin{equation}\label{eq:mM}
\boxed{\mathfrak m(r):=\sigma(r)\,k\,a\,(1+r^2)^{-\frac N2},\qquad
\mathfrak M(r):=\sigma(r)\,a\,\mathcal S(r),}
\end{equation}
that is,
\[
\mathfrak m=\min_{\partial\mathcal B_1}\sigma P_r \qquad\text{and}\qquad \mathfrak M=\max_{\partial\mathcal B_1}\sigma P_r.
\]
By Lemma~\ref{lem:extrema} the mean of $\sigma P_r$ on $\partial\mathcal B_1$ is $k\sigma(r)$; in particular
$\mathfrak m(r)<k\sigma(r)<\mathfrak M(r)$.

\subsection{The key estimate}
By Proposition~\ref{prop:normal}, the sign of $\partial_\nu u_\e$ at a point of $\partial\mathcal B_1\cap\{x_1=x_2=0\}$ and at a point $e_j$ is asymptotically the sign of $\mathfrak m(r_*)-1$ and of $\mathfrak M(r_*)-1$ respectively, where $r_*$ is the limit radius, which is a zero of $g$. Hence everything reduces to the following statement.

\begin{proposition}\label{prop:mM}
Let $N\ge5$ and $2\le k\le6$. If $g(r_*)=0$, then $\mathfrak m(r_*)<1<\mathfrak M(r_*)$.
\end{proposition}

Recall that, by \eqref{eq:g} and \eqref{eq:mw},
\begin{equation}\label{eq:gsigma}
g=1+k\beta\sigma-\varphi_k\,\sigma^{2-\ell}-\beta\,\sigma^{1-\ell} .
\end{equation}

\begin{lemma}\label{lem:sigma}
Let $N\ge5$, $k\ge2$ and $r\in(0,1)$. Then
\begin{enumerate}[label=\textup{(\roman*)}]
\item $0<\sigma<\sigma_0:=r^{-N}a^{N-1}$ and $\varphi_k<a^{-(N-2)}$;
\item $\sigma^{-1}\le a\,r^N\mathcal S(t)+\frac{\Lambda_k}2$, where $\Lambda_k$ is defined in Lemma~\ref{lem:zero};
\item if $2\le k\le6$, then $\Lambda_k\le k-1$ and $\sigma^{-1}\le r^Na^{-(N-1)}+(k-1)\,a\,r^N(a^2+r^2)^{-\frac N2}+\frac{k-1}2$.
\end{enumerate}
\end{lemma}

\begin{proof}
By \eqref{eq:s} and \eqref{eq:phiprime},
\begin{equation}\label{eq:sigmainv}
\frac1\sigma=\frac{r^{N-1}\varphi_k'}{2(N-2)}=r^Na^{-(N-1)}+\frac{r^{N-1}}{2(N-2)}\sum_{j=1}^{k-1}\omega_j ,
\end{equation}
and $\omega_j>0$, so $\sigma<\sigma_0$; moreover $\varphi_k<a^{-(N-2)}$ by \eqref{eq:phiexpl}, since $\tau_j>0$. By \eqref{eq:omegaj}, $\omega_j\le\frac{N-2}rb_j^{-(N-2)}+2(N-2)ra\,q_j^{-N/2}$, and $r^{N-2}b_j^{-(N-2)}=(2\mathfrak s_j)^{-(N-2)}$, $r^Na^{-(N-1)}=ar^Nq_0^{-N/2}$; inserting these in \eqref{eq:sigmainv} gives (ii). If $k\le6$, then $2\mathfrak s_j\ge2\sin\frac\pi k\ge1$ for $1\le j\le k-1$, so $\Lambda_k\le k-1$, $b_j\ge r$ and $q_j\ge a^2+r^2$, and (iii) follows from (ii).
\end{proof}

\subsubsection{The maximum}

\begin{lemma}\label{lem:fourier}
Let $N\ge5$, $k\ge2$ and $r\in(0,1)$. Then $a\big[\mathcal S(r)-r^N\mathcal S(t)\big]>k$.
\end{lemma}

\begin{proof}
By \eqref{eq:Fcos} with $\theta=0$ and $r$ replaced by $\rho\in(0,1)$,
\begin{equation}\label{eq:sumexp}
\mathcal S(\rho)=k\sum_{n\ge0}\sum_{l\in\mathbb Z}c_nc_{n+|l|k}\,\rho^{2n+|l|k}.
\end{equation}
Applying \eqref{eq:sumexp} with $\rho=r$ and with $\rho=t=r^2$ we obtain
\[
\mathcal S(r)-r^N\mathcal S(t)=k\sum_{n\ge0}\sum_{l\in\mathbb Z}c_nc_{n+|l|k}\,r^{2n+|l|k}\Big(1-r^{N+2n+|l|k}\Big).
\]
All the terms are nonnegative and $1-r^{N+2n+|l|k}\ge1-r^N\ge1-r^2=a$. Keeping only the terms with $l=0$ we get
\[
a\big[\mathcal S(r)-r^N\mathcal S(t)\big]\ge k\,a^2\sum_{n\ge0}c_n^2\,r^{2n}.
\]
Finally $c_n=\prod_{i=0}^{n-1}\frac{N/2+i}{1+i}\ge1$, hence
$\sum_nc_n^2r^{2n}\ge\sum_nc_nr^{2n}=a^{-\frac N2}$, and the left-hand side is at least $k\,a^{2-\frac N2}>k$, because $N>4$ and $a\in(0,1)$.
\end{proof}

\begin{proposition}\label{prop:M}
Let $N\ge5$ and $k\ge2$ be such that $\Lambda_k\le2k$. Then $\mathfrak M(r)>1$ for every $r\in(0,1)$. In particular, this holds for $2\le k\le6$.
\end{proposition}

\begin{proof}
By Lemma~\ref{lem:fourier} and Lemma~\ref{lem:sigma}\,(ii), $a\,\mathcal S(r)>a\,r^N\mathcal S(t)+k\ge a\,r^N\mathcal S(t)+\frac{\Lambda_k}2\ge\sigma^{-1}$, that is, $\mathfrak M(r)=\sigma a\mathcal S(r)>1$. If $k\le6$, then $\Lambda_k\le k-1<2k$ by Lemma~\ref{lem:sigma}\,(iii).
\end{proof}

\subsubsection{The minimum}
The inequality $\mathfrak m<1$ is false for $r$ small, so one has to show that the zeros of $g$ do not lie there. We prove the slightly stronger implication
\begin{equation}\label{eq:key}
\mathfrak m(r)\ge1\quad\Longrightarrow\quad g(r)>1\qquad(N\ge5,\ 2\le k\le6).
\end{equation}
Let
\[
\mathfrak t(t):=\frac{(1-t)^2}{t(1+t)},\qquad
\sigma_{\mathfrak m}(t):=\frac{(1+t)^{N/2}}{k(1-t)},\qquad\text{so that}\qquad \mathfrak m=\frac{\sigma}{\sigma_{\mathfrak m}(t)} ;
\]
$\mathfrak t$ is decreasing from $+\infty$ to $0$ on $(0,1)$, and $\sigma_{\mathfrak m}$ is increasing. By Lemma~\ref{lem:sigma}\,(i), $\sigma<\sigma_0=r^{-N}a^{N-1}=k\,\mathfrak t(t)^{N/2}\sigma_{\mathfrak m}(t)$. Hence, if $\mathfrak m(r)\ge1$,
\begin{equation}\label{eq:constraints}
\sigma\ \ge\ \sigma_{\mathfrak m}(t),\qquad k\,\mathfrak t(t)^{N/2}>1 .
\end{equation}
We denote by $\hat t_{N,k}$ the unique solution of $k\,\mathfrak t(t)^{N/2}=1$.

\begin{lemma}\label{lem:mmono}
Let $N\ge5$ and $k\ge2$. Then
\begin{equation}\label{eq:sigmaexact}
\frac1{\sigma(r)}=\frac{\Lambda_k}2+t^{\frac N2}B(t),\qquad
B(t):=\sum_{j=0}^{k-1}\frac{\cos\theta_j-t}{|1-te^{i\theta_j}|^{N}},
\end{equation}
where $B$ is nonnegative and increasing on $[0,1)$. Consequently $\mathfrak m$ is strictly decreasing on $(0,1)$, with $\mathfrak m(0^+)=\frac{2k}{\Lambda_k}$ and $\mathfrak m(1^-)=0$. If $\Lambda_k<2k$ (in particular if $2\le k\le6$), there is a unique $t_{N,k}\in(0,1)$ such that
\[
\mathfrak m(r)\ge1\iff r^2\le t_{N,k},\qquad\text{and}\qquad t_{N,k}<\hat t_{N,k}.
\]
\end{lemma}

\begin{proof}
By \eqref{eq:omegaj}, since $r^{N-2}b_j^{-(N-2)}=(2\mathfrak s_j)^{-(N-2)}$, $r^{N-2}b_j^2=4\mathfrak s_j^2r^N$ and $a-2\mathfrak s_j^2=\cos\theta_j-t$,
\[
\frac{r^{N-1}}{2(N-2)}\,\omega_j=\frac12\,(2\mathfrak s_j)^{-(N-2)}+t^{\frac N2}\,q_j^{-\frac N2}\,(\cos\theta_j-t).
\]
Inserting this in \eqref{eq:sigmainv}, and writing $r^Na^{-(N-1)}=t^{\frac N2}(1-t)\,q_0^{-\frac N2}$, we get \eqref{eq:sigmaexact}. Since $\frac{d}{dt}|1-te^{i\theta}|^2=-2(\cos\theta-t)$, we have
\[
B(t)=\frac1{N-2}\,\frac{d}{dt}\sum_{j=0}^{k-1}|1-te^{i\theta_j}|^{-(N-2)}
=\frac k{N-2}\sum_{n\ge0}\sum_{l\in\mathbb Z}(2n+|l|k)\,d_nd_{n+|l|k}\,t^{2n+|l|k-1},
\]
where $d_n:=\frac{((N-2)/2)_n}{n!}>0$ and the last identity is \eqref{eq:sumexp} with $N$ replaced by $N-2$. Hence $B\ge0$ and $B$ is increasing, so $1/\sigma$ is increasing. Then
\[
\frac1{\mathfrak m}=\frac{(1+t)^{\frac N2}}{k(1-t)}\cdot\frac1\sigma
\]
is the product of two positive increasing functions, the first one strictly increasing. Finally $1/\sigma\to\Lambda_k/2$ as $t\to0$, and $\mathfrak m<k\,\mathfrak t(t)^{N/2}\to0$ as $t\to1$. If $\Lambda_k<2k$ then $\mathfrak m(0^+)>1$, and $t_{N,k}$ exists and is unique; moreover $1=\mathfrak m(\sqrt{t_{N,k}})<k\,\mathfrak t(t_{N,k})^{N/2}$ gives $t_{N,k}<\hat t_{N,k}$. For $2\le k\le6$ we have $\Lambda_k\le k-1$ by Lemma~\ref{lem:sigma}\,(iii).
\end{proof}

\noindent For $N\ge6$, with $m$ and $\ell$ as in \eqref{eq:mw}, we set
\[
T_{N,k}(t):=\big(1-t^{m}\big)\Big[k-\Big(\frac{k(1-t)}{1+t}\Big)^\ell(1+t)^{-2}\Big]-\frac{t}{(1-t)^3}\Big(\frac t{1-t}\Big)^\ell .
\]

\begin{lemma}\label{lem:Tconc}
Let $N\ge6$. As a function of $t\in[0,\frac12]$ and of the real parameter $k\ge1$, $T_{N,k}$ is concave in $t$, increasing in $k$, and $T_{N,k}(0)=k-k^\ell\ge0$.
\end{lemma}

\begin{proof}
Clearly $T_{N, k}(0)=k-k^\ell\geq 0$, since $k\geq 1$ and $\ell\le1$. Moreover
$$\partial_k T_{N, k}(t)=(1-t^m)\left[1-\ell k^{\ell-1}\frac{(1-t)^\ell}{(1+t)^{\ell+2}}\right]>0$$ for $t\in(0,\frac12]$, since $k\geq 1$, $\ell\le1$ (so that $\ell k^{\ell-1}\le1$) and $\frac{(1-t)^\ell}{(1+t)^{\ell+2}}<1$. Finally, a direct computation shows that \begin{equation}\label{eq:T2}
\begin{aligned}
T_{N,k}''(t)={}&-m(m-1)\,t^{m-2}\Big[k-k^\ell(1-t)^\ell(1+t)^{-\ell-2}\Big]\\
&-4m\,k^\ell\,t^{m-1}(\ell+1-t)\,(1-t)^{\ell-1}(1+t)^{-\ell-3}\\
&-\big(1-t^m\big)\,k^\ell\,(1-t)^{\ell-2}(1+t)^{-\ell-4}\Big[6t^2-12(\ell+1)t+4\ell^2+8\ell+6\Big]\\
&-t^{\ell-1}(1-t)^{-\ell-5}\Big[\ell(\ell+1)+6(\ell+1)t+6t^2\Big].
\end{aligned}
\end{equation}
All the four terms in \eqref{eq:T2} are nonpositive for $t\in(0,\frac12]$: indeed $k-k^\ell(1-t)^\ell(1+t)^{-\ell-2}\ge k-k^\ell\ge0$, $\ell+1-t>0$, and the polynomial $6t^2-12(\ell+1)t+4\ell^2+8\ell+6$ is decreasing for $t<\ell+1$, so that on $[0,\frac12]$ it is bounded from below by its value at $t=\frac12$, namely $4\ell^2+2\ell+\frac32>0$. Hence $T_{N,k}$ is concave on $(0,\frac12]$.
\end{proof}
\begin{proposition}\label{prop:key6}
Let $N\ge6$ and $2\le k\le6$. Then \eqref{eq:key} holds.
\end{proposition}

\begin{proof}
Let $\ell=\frac4{N-2}\in(0,1]$ and $m=\frac{N-2}2$ be as in \eqref{eq:mw}, and assume $\mathfrak m(r)\ge1$.

\smallskip\noindent\emph{Step 1: a lower bound for $g$.} By \eqref{eq:gsigma},
\[
g=1+\sigma\big[k\beta-\varphi_k\,\sigma^{1-\ell}-\beta\,\sigma^{-\ell}\big].
\]
We estimate from above the two negative terms in the bracket. Since $-\ell<0$, the first inequality in \eqref{eq:constraints} gives
\[
\sigma^{-\ell}\le\sigma_{\mathfrak m}(t)^{-\ell}=\Big(\frac{k(1-t)}{1+t}\Big)^\ell(1+t)^{-2},
\]
where we used $\frac{N\ell}{2}=\ell+2$. Since $1-\ell\ge0$, Lemma~\ref{lem:sigma}\,(i) gives $\sigma^{1-\ell}\le\sigma_0^{1-\ell}$; hence, if $\varphi_k>0$, using also $\varphi_k<a^{-(N-2)}$,
\[
\varphi_k\,\sigma^{1-\ell}\le a^{-(N-2)}\sigma_0^{1-\ell},
\]
and this inequality is trivially true if $\varphi_k\le0$, because its right-hand side is positive. Finally $\beta=r^{-(N-2)}(1-t^m)$ and, since $\sigma_0=r^{-N}a^{N-1}$ and $N\ell=4+2\ell$,
\[
a^{-(N-2)}\sigma_0^{1-\ell}=r^{-(N-2)}\,\frac t{(1-t)^3}\Big(\frac t{1-t}\Big)^\ell .
\]
Inserting these estimates in the bracket we obtain
\begin{equation}\label{eq:Tbound}
g\ \ge\ 1+\sigma\,r^{-(N-2)}\,T_{N,k}(t),
\end{equation}
so that it suffices to show that $T_{N,k}(t)>0$.

\smallskip\noindent\emph{Step 2: localisation.} Let
\[
\bar t:=\begin{cases}\frac9{25} & \text{if }k=2,\\[2pt] \frac25 & \text{if }3\le k\le6,\end{cases}\quad\text{when }N\ge11,
\qquad
\bar t:=\begin{cases}\frac38 & \text{if }k=2,\\[2pt] \frac{17}{40} & \text{if }3\le k\le6,\end{cases}\quad\text{when }6\le N\le10 .
\]
We claim that $k\,\mathfrak t(\bar t)^{N/2}<1$. Since $\mathfrak t(\bar t)<1$, the left-hand side is increasing in $k$ and decreasing in $N$, so it suffices to consider $k=2$, $k=6$ and the smallest $N$ in each range; as $\mathfrak t(\frac9{25})=\frac{128}{153}$, $\mathfrak t(\frac25)=\frac9{14}$, $\mathfrak t(\frac38)=\frac{25}{33}$ and $\mathfrak t(\frac{17}{40})=\frac{529}{969}$, a direct computation gives that in every case $k\,\mathfrak t(\bar t)^{N/2}<1$.
Since $\mathfrak t$ is decreasing, this gives $\hat t_{N,k}<\bar t$. By Lemma~\ref{lem:mmono}, the assumption $\mathfrak m(r)\ge1$ is equivalent to $t\le t_{N,k}$, and $t_{N,k}<\hat t_{N,k}$. Therefore
\[
0<t\le t_{N,k}<\hat t_{N,k}<\bar t<\tfrac12 .
\]

\smallskip\noindent\emph{Step 3: reduction to the endpoint.} By Lemma~\ref{lem:Tconc}, $T_{N,k}$ is concave on $[0,\frac12]$ and $T_{N,k}(0)=k-k^\ell\ge0$. Hence, for $0<t\le\bar t$, $T_{N,k}$ lies above the chord joining $(0,T_{N,k}(0))$ and $(\bar t,T_{N,k}(\bar t))$:
\[
T_{N,k}(t)\ \ge\ \Big(1-\frac t{\bar t}\Big)T_{N,k}(0)+\frac t{\bar t}\,T_{N,k}(\bar t)\ \ge\ \frac t{\bar t}\,T_{N,k}(\bar t).
\]
Since $T_{N,k}$ is increasing in $k$ (Lemma~\ref{lem:Tconc}) and $\bar t$ does not depend on $k\in\{3,\dots,6\}$, we have $T_{N,k}(\bar t)\ge T_{N,k_0}(\bar t)$, with $k_0=2$ if $k=2$ and $k_0=3$ if $3\le k\le6$. Finally, a direct computation, in which $1-\bar t^{\,m}$ is bounded from below by its value at the smallest $m$ of each range and the powers of exponent $\ell$ by their values at the endpoints of the range of $\ell$, gives
\[
T_{N,2}\big(\tfrac9{25}\big)>0.06,\qquad T_{N,3}\big(\tfrac25\big)>0.5\qquad(N\ge11),
\]
\[
T_{N,2}\big(\tfrac38\big)>0.07,\qquad T_{N,3}\big(\tfrac{17}{40}\big)>0.04\qquad(6\le N\le10).
\]
Hence $T_{N,k}(t)>0$, and $g(r)>1$ by \eqref{eq:Tbound}.
\end{proof}

\begin{proposition}\label{prop:key5}
Let $N=5$ and $2\le k\le6$. Then \eqref{eq:key} holds.
\end{proposition}

\begin{proof}
For $N=5$ we have $\ell=\frac43$, so that $2-\ell=\frac23$ and $1-\ell=-\frac13$. Let
\[
h_k(v):=k\,v^{1/3}-v^{-1},\qquad \Upsilon(r):=\frac{a^{-3}}{\beta}=\frac{r^3}{(1-r^3)(1-r^2)^3},
\]
and
\[
T_{5,k}(t):=(3k^2)^{\frac13}-\frac k3-\Upsilon\big(\sqrt t\,\big)=(3k^2)^{\frac13}-\frac k3-\frac{t^{3/2}}{(1-t^{3/2})(1-t)^3}.
\]
The functions $h_k$ and $\Upsilon$ are increasing on $(0,\infty)$ and on $(0,1)$ respectively; in particular $T_{5,k}$ is decreasing on $(0,1)$. Assume $\mathfrak m(r)\ge1$.

\smallskip\noindent\emph{Step 1: a lower bound for $g$.} By \eqref{eq:gsigma},
\[
g=1+k\beta\sigma-\varphi_k\sigma^{\frac23}-\beta\sigma^{-\frac13}=1+\sigma^{\frac23}\big[\beta\,h_k(\sigma)-\varphi_k\big],
\]
and by Lemma~\ref{lem:sigma}\,(i), $\varphi_k<a^{-3}=\beta\,\Upsilon(r)$ (this is trivial if $\varphi_k\le0$). Since $1-\ell<0$, here both negative terms of $g$ have to be estimated with a lower bound for $\sigma$. We first show that
\begin{equation}\label{eq:m13}
\mathfrak m\big(3^{-1/2}\big)\ge1 .
\end{equation}
At $t=\frac13$ we have $\mathfrak m=\frac{3\sqrt3}{16}\,k\,\sigma$. In \eqref{eq:sigmaexact} the terms with $\cos\theta_j\le\frac13$ are nonpositive, and only $j=0$ and, for $k=6$, $j=1,5$ have $\cos\theta_j>\frac13$; dropping the other terms, a direct computation gives
\[
\frac1{\sigma}\le\frac{\Lambda_k}2+\frac{3\sqrt3}{16}\,(1+\delta_k),\qquad \delta_k=0\ \ (k\le5),\qquad \delta_6=\tfrac{16}{243}\big(\tfrac97\big)^{5/2}<\tfrac18 ,
\]
and \eqref{eq:m13} follows from $\Lambda_k<\frac{16}{25}(k-1-\delta_k)$, which is true since $\Lambda_2=\frac18$, $\Lambda_3<\frac25$, $\Lambda_4<\frac56$, $\Lambda_5<\frac{31}{20}$ and $\Lambda_6<\frac{101}{40}$.

By Lemma~\ref{lem:mmono}, the assumption $\mathfrak m(r)\ge1$ is equivalent to $t\le t_{5,k}$, and \eqref{eq:m13} means $\frac13\le t_{5,k}$. Since $1/\sigma$ is increasing, $\mathfrak m=\sigma/\sigma_{\mathfrak m}$ equals $1$ at $t_{5,k}$, and $\sigma_{\mathfrak m}$ is increasing, we get
\[
\sigma(r)\ \ge\ \sigma\big(\sqrt{t_{5,k}}\big)=\sigma_{\mathfrak m}(t_{5,k})\ \ge\ \sigma_{\mathfrak m}\big(\tfrac13\big)=\frac{16\sqrt3}{9k}\ >\ \frac3k .
\]
As $h_k$ is increasing, $h_k(\sigma(r))\ge h_k(\frac3k)=(3k^2)^{1/3}-\frac k3$, and therefore
\begin{equation}\label{eq:T5bound}
g\ \ge\ 1+\sigma^{\frac23}\,\beta\,\big[h_k(\sigma)-\Upsilon(r)\big]\ \ge\ 1+\sigma^{\frac23}\,\beta\,T_{5,k}(t),
\end{equation}
so that it suffices to show that $T_{5,k}(t)>0$.

\smallskip\noindent\emph{Step 2: localisation.} Let
\[
(\bar t_2,\bar t_3,\bar t_4,\bar t_5,\bar t_6):=\Big(\tfrac{19}{50},\tfrac25,\tfrac{21}{50},\tfrac{87}{200},\tfrac{89}{200}\Big).
\]
A direct computation gives $k^2\,\mathfrak t(\bar t_k)^5<1$, i.e.\ $k\,\mathfrak t(\bar t_k)^{5/2}<1$, so that $\hat t_{5,k}<\bar t_k$, because $\mathfrak t$ is decreasing. By Lemma~\ref{lem:mmono}, $t\le t_{5,k}<\hat t_{5,k}$. Therefore
\[
0<t\le t_{5,k}<\hat t_{5,k}<\bar t_k .
\]

\smallskip\noindent\emph{Step 3: reduction to the endpoint.} Since $T_{5,k}$ is decreasing, $T_{5,k}(t)\ge T_{5,k}(\bar t_k)$ for $0<t\le\bar t_k$, and a direct computation gives
\[
T_{5,k}(\bar t_k)>0.29\qquad\text{for }2\le k\le6 .
\]
Hence $T_{5,k}(t)>0$, and $g(r)>1$ by \eqref{eq:T5bound}.
\end{proof}

\begin{proof}[Proof of Proposition~\ref{prop:mM}]
$\mathfrak M(r_*)>1$ by Proposition~\ref{prop:M}. If $\mathfrak m(r_*)\ge1$, then $g(r_*)>1$ by Propositions~\ref{prop:key6} and~\ref{prop:key5}, which contradicts $g(r_*)=0$.
\end{proof}

\subsection{Proof of Theorem~\ref{thm:nodal}}

Let $u_\e=V_{y_\e}+\phi_{y_\e}$ be as in Theorem~\ref{thm:nodal}. In Corollary~\ref{cor:exist} the parameters $y_\e$ lie in a fixed compact set $\overline{\mathfrak U_c}\subset\interior\mathcal P_{\eta'}$, and every limit point $y_*=(r_*,d_{0*},d_{1*})$ of $y_\e$ is a critical point of $\Psi_{A,B}$, so that $g(r_*)=0$. By Proposition~\ref{prop:normal}, along a sequence $\e_n\to0$ with $y_{\e_n}\to y_*$,
\begin{equation}\label{eq:normlim}
\mu_n^{-\frac{N-2}{2}}\,\partial_\nu u_{\e_n}\to (N-2)\,\alpha_N\,d_{0*}^{\frac{N-2}{2}}\big(\sigma(r_*)P_{r_*}-1\big)\qquad\text{uniformly on }\partial\mathcal B_1 .
\end{equation}

\begin{proof}[Proof of Theorem~\ref{thm:nodal}]
(a) Otherwise there are $\e_n\to0$ and $x_n\in\partial\mathcal B_1\cap\{x_1=x_2=0\}$ with
$\partial_\nu u_{\e_n}(x_n)\ge0$. Up to a subsequence, $y_{\e_n}\to y_*$ and, by \eqref{eq:normlim}, Lemma~\ref{lem:extrema} and Proposition~\ref{prop:mM}, $\mu_n^{-\frac{N-2}{2}}\partial_\nu u_{\e_n}(x_n)\to(N-2)\alpha_Nd_{0*}^{\frac{N-2}{2}}\big(\mathfrak m(r_*)-1\big)<0$, a contradiction. Since $u_\e=0$ on $\partial\mathcal B_1$ and
$\partial_\nu u_\e<0$ on the compact set $\partial\mathcal B_1\cap\{x_1=x_2=0\}$, $u_\e>0$ in a one-sided
neighbourhood of it.

(b) The same argument, using $\sigma(r_*)P_{r_*}(e_j)=\mathfrak M(r_*)>1$.

(c) Let $\gamma:[0,1]\to\partial\mathcal B_1$ be a path from a point of $\partial\mathcal B_1\cap\{x_1=x_2=0\}$ to $e_0$. By
(a), (b) and the continuity of $\nabla u_\e$ up to $\partial\mathcal B_1$, for $\lambda<1$ close to $1$ we have
$u_\e(\lambda\gamma(0))>0>u_\e(\lambda\gamma(1))$, so $u_\e$ vanishes at some point
of $\lambda\gamma([0,1])$, at distance $1-\lambda$ from $\partial\mathcal B_1$. Letting $\lambda\to1$ we get a
point of $\mathcal N(u_\e)\cap\partial\mathcal B_1$.
\end{proof}

\begin{remark}
More precisely, \eqref{eq:normlim} shows that every limit point, as $\e\to0$, of zeros of $\partial_\nu u_\e$ on $\partial\mathcal B_1$ belongs to the level set $\{x\in\partial\mathcal B_1:\ \sigma(r_*)P_{r_*}(x)=1\}$, which separates $\partial\mathcal B_1\cap\{x_1=x_2=0\}$ from the points $e_j$.
\end{remark}

\begin{remark}\label{rem:twosol}
The argument above also describes how the nodal set meets the plane $\mathcal X=\{x_3=\dots=x_N=0\}$ of the peaks. Let $F$ and $\mathbb D$ be as in the proof of Lemma~\ref{lem:extrema}, so that $P_r(x)=a\,F(z)$ for $x=(z,x'')\in\partial\mathcal B_1$, and set
\[
\mathfrak m_{\mathcal X}(r):=\min_{\partial\mathcal B_1\cap\mathcal X}\sigma(r)P_r=\sigma(r)\,a\min_{\theta\in\R}F(e^{i\theta}),
\]
so that $\mathfrak m(r)\le\mathfrak m_{\mathcal X}(r)<\mathfrak M(r)$, with $\mathfrak m_{\mathcal X}=\mathfrak m$ when $k=2$. If $\mathfrak m_{\mathcal X}(r_*)>1$, then by \eqref{eq:normlim} $\partial_\nu u_\e>0$ on the circle $\partial\mathcal B_1\cap\mathcal X$ for $\e$ small, so that $u_\e<0$ near it; since $u_\e>0$ near the origin, in the plane $\mathcal X$ the nodal set separates the origin from $\partial\mathcal B_1\cap\mathcal X$ and does not meet it. If instead $\mathfrak m_{\mathcal X}(r_*)<1$, the argument used in the proof of Theorem~\ref{thm:nodal}\,(c) shows that the nodal set meets $\partial\mathcal B_1\cap\mathcal X$. Moreover, since $F$ is convex, the set $\mathcal C:=\{z\in\overline{\mathbb D}:\ \sigma(r_*)\,a\,F(z)\le1\}$ is convex, and it contains a neighbourhood of $0$ by Proposition~\ref{prop:mM}. If $\mathfrak m_{\mathcal X}(r_*)>1$, then $\mathcal C\subset\mathbb D$, and the level set $\{x\in\partial\mathcal B_1:\ \sigma(r_*)P_{r_*}(x)=1\}=\{(z,x'')\in\partial\mathcal B_1:\ z\in\partial\mathcal C\}$ is homeomorphic to $\mathbb S^1\times\mathbb S^{N-3}$.

Numerically (we checked $6\le N\le40$ and $N=50,60,80,100$), for $3\le k\le6$ the minimum of $F$ on $\partial\mathbb D$ is attained at the midpoints $e^{i(\theta_j+\pi/k)}$ of the arcs between consecutive peaks, and, denoting by $r^-<r^+$ the two zeros of $g$:
\begin{enumerate}[label=\textup{(\roman*)}]
\item $\mathfrak m_{\mathcal X}(r^-)>1$ in all cases, with $\mathfrak m_{\mathcal X}(r^-)\ge1.87$; for $N=6$ we find $\mathfrak m_{\mathcal X}(r^-)\approx1.88$, $4.33$, $7.58$, $10.87$ for $k=3,4,5,6$;
\item $\mathfrak m_{\mathcal X}(r^+)<1$ if and only if $k\in\{3,4\}$, or $k=5$ and $N\le9$, or $k=6$ and $N\le7$; for $N=6$ we find $\mathfrak m_{\mathcal X}(r^+)\approx4\cdot10^{-4}$, $2\cdot10^{-3}$, $7\cdot10^{-3}$, $2\cdot10^{-2}$ for $k=3,4,5,6$, and the nodal set meets $\partial\mathcal B_1\cap\mathcal X$ at angular distance about $12^\circ$ from each point $e_j$.
\end{enumerate}
Also, for $N=5$ and $3\le k\le6$ the zero of $g$ satisfies $\mathfrak m_{\mathcal X}>1$. Finally, since $u_\e$ depends only on $(x_1,x_2,|x''|)$ and the orbits of $O(N-2)$ are connected, the nodal domains of $u_\e$ correspond to those of a function of three variables. Counting on a grid the positive and negative components of the limit profile $G(\cdot,0)-\sigma(r_*)\sum_jG(\cdot,\xi_j(r_*))$, we find one positive and one negative component at $r_*=r^-$ (and for $N=5$, $3\le k\le6$), and one positive and $k$ negative components at $r_*=r^+$ in the range described in (ii); for $k=2$ we find one positive and two negative components at both zeros. This is the origin of the description given in Remark~\ref{rem:two}.
\end{remark}

\subsection{The number of nodal domains}
Theorem~\ref{thm:nodal} concerns the intersection of the nodal set with the boundary. The number of nodal domains of $u_\e$, i.e.\ of connected components of $\mathcal B_1\setminus\{u_\e=0\}$, can be partially controlled by a simple energy argument.

\begin{proposition}\label{prop:domains}
Let $N\ge5$, $k\ge2$, and let $u_\e$ be a family of solutions of \eqref{pb} of the form \eqref{eq:form} given by Theorem~\ref{thm:main1} or Theorem~\ref{thm:main2}\,\textup{(i)}. Then, for $\e$ small enough, $u_\e$ has exactly one positive nodal domain and at most $k$ negative nodal domains; in particular, it has at least $2$ and at most $k+1$ nodal domains.
\end{proposition}

\begin{proof}
Let $\mathcal D$ be a nodal domain of $u_\e$ and $\psi:=u_\e\mathbf 1_{\mathcal D}\in H^1_0(\mathcal B_1)$. Testing \eqref{pb} with $\psi$ and using the Poincar\'e inequality ($\e<\lambda_1(\mathcal B_1)$) and the Sobolev inequality $S\|\psi\|_{L^{2^*}}^2\le\|\nabla \psi\|_{L^2}^2$, we get
\[
\Big(1-\frac{\e}{\lambda_1}\Big)\int|\nabla \psi|^2\le\int|\psi|^{2^*}\le\Big(S^{-1}\int|\nabla \psi|^2\Big)^{\frac N{N-2}},
\]
so that
\[
\int_{\mathcal D}|\nabla u_\e|^2\ge\Big(1-\frac{\e}{\lambda_1}\Big)^{\frac{N-2}2}S^{\frac N2}.
\]
On the other hand, testing \eqref{pb} with $u_\e^\pm$,
\[
\int|\nabla u_\e^\pm|^2=\int|u_\e^\pm|^{2^*}+\e\int|u_\e^\pm|^2 .
\]
Write $u_\e=V_\e+\phi_\e$, with $V_\e=PU_{\delta_0,0}-\sum_jPU_{\delta_1,\xi_j}$ and $\|\phi_\e\|_{H^1_0}\to0$. Since $\zeta\mapsto\zeta^\pm$ is $1$-Lipschitz, $\|u_\e^\pm-V_\e^\pm\|_{L^{2^*}}\le\|\phi_\e\|_{L^{2^*}}\to0$. Moreover $0\le PU_{\delta,\xi}\le U_{\delta,\xi}$, hence $V_\e^+\le U_{\delta_0,0}$ and $V_\e^-\le\sum_jU_{\delta_1,\xi_j}$, so that
\[
\int|V_\e^+|^{2^*}\le\int_{\R^N}U^{2^*}=S^{\frac N2},\qquad
\int|V_\e^-|^{2^*}\le\int_{\R^N}\Big(\sum_jU_{\delta_1,\xi_j}\Big)^{2^*}=k\,S^{\frac N2}+o(1),
\]
because the concentration points are at mutual distance bounded away from zero. Therefore $\int|\nabla u_\e^+|^2\le S^{N/2}+o(1)$ and $\int|\nabla u_\e^-|^2\le kS^{N/2}+o(1)$. Since the energy of $u_\e^\pm$ is the sum of the energies of the positive, resp.\ negative, nodal domains, each of which is at least $(1-o(1))S^{N/2}$, there are at most one positive and at most $k$ negative nodal domains for $\e$ small. Finally, $u_\e^\pm\not\equiv0$ for $\e$ small. Indeed, $V_\e$ coincides with $U_{\delta_0,0}$ (resp.\ with $-U_{\delta_1,\xi_j}$), up to a term which is $o(1)$ in $L^{2^*}$, in a fixed ball centred at $0$ (resp.\ at $\xi_j$), so that
\[
\int|V_\e^+|^{2^*}\ge S^{\frac N2}+o(1),\qquad \int|V_\e^-|^{2^*}\ge k\,S^{\frac N2}+o(1),
\]
and $\|u_\e^\pm-V_\e^\pm\|_{L^{2^*}}\to0$. Hence $u_\e$ has at least one positive and at least one negative nodal domain.
\end{proof}

\begin{remark}[On the restriction $k\le6$]\label{rem:k6}
The inequality $\mathfrak M>1$ holds whenever $\Lambda_k\le2k$ (Proposition~\ref{prop:M}), which is true for some $k\ge7$ as well. The bound $k\le6$ is really used in the proof of $\mathfrak m<1$ (Propositions \ref{prop:key6} and \ref{prop:key5}).
Note that $\mathfrak M>1$ genuinely fails for $k$ large and $r$ small: there
$\sigma\to2/\Lambda_k$ and $\mathfrak M\to2k/\Lambda_k\to0$ as $k\to\infty$. Hence, for $k$ large, one also needs
to localise the zeros of $g$ away from $r=0$. For $N=5$ this is provided by Proposition~\ref{prop:N5ring}, which, however, is not quantitative enough to cover small values of $k$.
\end{remark}

\section*{Data Availability Statements}
All data generated or analysed during this study are included in this article.
\section*{Declarations}
\textbf{Conflicts of Interest:} The authors declare they have no financial interests.
\section*{Acknowledgements}
The authors acknowledge the use of ChatGPT (Open AI) to assist selected calculations and
the drafting and revision of the manuscript. They have independently verified all mathematical arguments and calculations to take full responsibility for the content of the paper.

\end{document}